\documentclass{gOMS2e}

\usepackage{amssymb,amsmath}
\usepackage{graphicx,booktabs}
\usepackage{color}
\def\RR{\mathbb{R}}

\theoremstyle{plain}
\newtheorem{theorem}{Theorem}[section]

\newtheorem{lemma}[theorem]{Lemma}

\theoremstyle{definition}

\theoremstyle{algorithm}
\newtheorem{algorithm}{Algorithm}

\theoremstyle{remark}
\newtheorem{remark}{Remark}
\newenvironment{Proof}{\begin{trivlist}\item[]
  {\em Proof.\ }}{$\square$\end{trivlist}}

\date{}

\begin{document}

\title{A first-order multigrid method for bound-constrained\\ convex optimization}
\author{
\name{Michal Ko\v{c}vara$^{\rm a}$$^{\ast}$\thanks{$^\ast$Corresponding author. Email: m.kocvara@bham.ac.uk} and
Sudaba Mohammed$^{\rm b}$}
\affil{$^{a}$School of Mathematics, University of
    Birmingham, Edgbaston, Birmingham B15 2TT, UK, and Institute of Information Theory
and Automation, Czech Academy of Sciences, Pod
vod\'arenskou v\v{e}\v{z}\'{\i}~4, 18208 Praha 8, Czech Republic\\
  $^{b}$School of Mathematics, University of Birmingham, Edgbaston,
    Birmingham B15 2TT, UK,\\ on leave from Department of Mathematics, University of Kirkuk, Iraq}
}
\maketitle
\begin{abstract}
The aim of this paper is to design an efficient multigrid method for constrained convex optimization problems arising from discretization of some underlying infinite dimensional problems. Due to problem dependency of this approach, we only consider bound constraints with (possibly) a single equality constraint. As our aim is to target large-scale problems, we want to avoid computation of second derivatives of the objective function, thus excluding Newton like methods. We propose a smoothing operator that only uses first-order information and study the computational efficiency of the resulting method.
\end{abstract}


\begin{classcode}90C30, 65N55, 90C33\end{classcode}
\section{Introduction}
Multigrid methods were originally developed for the solution of large systems of linear algebraic equations arising from discretization of partial differential equations. Their practical utility and efficiency were demonstrated by Achi Brandt in his pioneering papers
\cite{brandt1973multi, brandt1977multi}. It has been well-known since the dark ages of multigrid that the methods may lose their superior efficiency when used for slightly different type of problems, namely the linear complementarity problems (LCP) \cite{brandt1983multigrid,hackbusch1983,mandel1984multilevel}. This is caused by the presence of unilateral obstacles (or box constraints in the optimization formulation of the problem). The fact that the sets of active constraints may vary for different discretization levels, and that the constraints may not even be recognized on very coarse meshes, may lead to poor quality of the coarse level corrections and, in effect, to significant deterioration or even loss of convergence of the method. Various remedies have been proposed by different authors \cite{brandt1983multigrid,hackbusch1983,hoppe1987multigrid,hoppe1987two,mandel1984multilevel}; these usually resulted in ``conservative" methods that were often significantly slower than standard methods for linear systems. Finally Kornhuber \cite{Kornhuber1994} proposed a truncated monotone multigrid method for LCP problems. This method has the property that as soon as the set of active constraints is correctly identified, the method converges with the same speed as without the presence of the constraints.

Not many attempts have been done to generalize the multigrid technique to the solution of optimization problems. From the successful ones, most focused on unconstrained problems \cite{frandi2013coordinate,frandi2015,gratton2010numerical,gratton2008recursive,lewis2000multigrid,lewis2013multigrid,nash2000multigrid}. In this case, the problem can be often identified with a discretized nonlinear PDE and thus techniques of nonlinear multigrid can be used.
These techniques use almost exclusively a variant of the Newton method or the Newton-Gauss-Seidel method as a smoother. Hence they require computation and storage of the Hessian of the underlying optimization problem at each iteration. This may be very costly or even prohibitive for some large-scale problems. Moreover, one step of the Newton-Gauss-Seidel method has high computational complexity and cannot be directly parallelized. Therefore, our goal is to use a smoother that only relies on first-order information of the optimization problem. Used only on the finest grid, it may be very inefficient, as compared to the Newton method, but this is where the coarse grid correction will help, just as in the linear case.

Treating general (equality or inequality) constraints by multigrid may be difficult, if not impossible, as we may not be able to find the corresponding restriction operators. If the number of constraints is directly proportional to the number of variables (such as for the bound constraints), the restriction operator for these constraints may be based on that for the variables. On the other hand, if the number of constraints is independent of the discretization (e.g., a single equality constraint) then the prolongation/restriction is simply the identity. All other situations are difficult, in our opinion. For this reason, all articles on multigrid for constrained problems either treat the bound-constrained problems or problems with a single equality constraint (e.g., \cite{gelman-mandel,gratton2010numerical,gratton2008recursive_bound,vallejos2010mgopt}) or assume that a restriction operator for the constraints exists \cite{nash2010convergence}.

The first, straightforward goal of this paper is to extend Kornhuber's technique to nonlinear convex optimization problems with bound constraints. The authors are not aware of this generalization in the existing literature. The second goal is to propose a smoothing operator that would only use first-order information, and study the efficiency of the resulting method. Finally, we extend the developed algorithm to problems with an additional equality constraint. We study the behaviour of the proposed algorithms on a number of numerical examples.

The paper is structured as follows. In Section 2 we introduce three multigrid algorithms for bound-constrained convex optimization problems; first Kornhuber's truncated correction scheme multigrid for bound-constrained quadratic problems, then its generalization to convex problems using the full approximation scheme and finally a version of the latter algorithm without truncation. This version is a new though minor modification of an existing method. Section 3 shows that the third algorithm can be easily extended to problems with an additional linear equality constraint that includes all variables. In Section 4 we introduce the first-order smoother---the gradient projection method with a gradient-based line search---and analyse its smoothing properties. Section 5 is devoted to numerical experiments: we, in particular, focus on the dependence of the convergence rate on the number of refinement levels and the number of the smoothing steps.
\section{Multigrid for bound-constrained optimization}
\subsection{The problem}
We consider the bound-constrained nonlinear optimization problem
\begin{align}
  &\min_{x\in\RR^n} f(x)\label{eq:1}\\
  &\mbox{subject to}\nonumber\\
  &\qquad \varphi_i\leq x_i \leq \psi_i,\quad i=1,\ldots,n\,,\nonumber
\end{align}
where $f$ is a \emph{convex} continuously differentiable nonlinear function and $\varphi_i<\psi_i$ for all~$i$. To guarantee existence of a solution, we assume that either \emph{the function $f$ is coercive} or that \emph{the feasible set}
$$
  {\cal F}=\left\{x\in\RR^n \mid \varphi_i\leq x_i\leq\psi_i,\ i=1,\ldots,n\right\}
$$
\emph{is bounded}.

In order to define and use a multigrid method, we assume that there is
a nested sequence of optimization problems arising from finite element
discretizations of some underlying problem:
\begin{align}
  &\min_{x^{(k)}\in\RR^{n_k}} f_k(x^{(k)})\label{eq:2}\\
  &\mbox{subject to}\nonumber\\
  &\qquad \varphi_i^{(k)}\leq x_i^{(k)} \leq \psi_i^{(k)},\quad i\in{\cal I}_k\,,\nonumber
\end{align}
where $k=0,1,\ldots,j$. Here $j$ is the finest discretization level
corresponding to the original problem (\ref{eq:1}) and 0 is the
coarsest level. The set ${\cal I}_k$ contains indices of $n_k$ finite
element vertices on discretization level $k$ and $f_k$ are the
discretizations of the same infinite dimensional nonlinear function.

The communication between two consecutive levels $k-1$ and $k$ is maintained by the prolongation and restriction operators $I_{k-1}^k:\RR^{n_{k-1}}\to\RR^{n_k}$ and $I_{k}^{k-1}:\RR^{n_{k}}\to\RR^{n_{k-1}}$, respectively, whereas $I_{k}^{k-1}=c(I_{k-1}^k)^T$. For the specific choice of these operators, see the last section.

For a feasible vector $x^{(k)}$ we introduce the sets of active
indices:
\begin{align*}
{\cal A}^\varphi_k(x^{(k)}) &:= \left\{i\in {\cal I}_k \mid x_i^{(k)} = \varphi_i^{(k)}\right\}\\
{\cal A}^\psi_k(x^{(k)}) &:= \left\{i\in {\cal I}_k \mid x_i^{(k)} = \psi_i^{(k)}\right\}\,.
\end{align*}
We will always assume that
$$
  {\cal A}^\varphi_k(x^{(k)}) \cap {\cal A}^\psi_k(x^{(k)}) = \emptyset\,.
$$
Finally, we denote the set of free variables by
$$
  {\cal A}^0_k := {\cal I}_k\setminus({\cal A}^\varphi_k(x^{(k)}) \cup {\cal A}^\psi_k(x^{(k)}))\,.
$$

Let us now introduce a \emph{modified} problem
\begin{align}
  &\min_{x^{(k)}\in\RR^{n_k}} f_k(x^{(k)})-v_k^T x^{(k)}\label{eq:2a}\\
  &\mbox{subject to}\nonumber\\
  &\qquad \varphi_i^{(k)}\leq x_i^{(k)} \leq \psi_i^{(k)},\quad i\in{\cal I}_k\,,\nonumber
\end{align}
where $v_k\in\RR^{n_k}$, together with the operator
$$
  x^{(k)}_{\rm new} = opt(f_k,\varphi^{(k)},\psi^{(k)},v_k;x^{(k)},\epsilon,\nu)
$$
that will play a role of the smoother for (\ref{eq:2a}).
The input are the problem data and the initial point $x^{(k)}$, the required
precision $\epsilon$ and the maximum number of iterations allowed $\nu$. We assume that $opt$ is a descent convergent algorithm for
the solution of problem (\ref{eq:2a}). In Section~\ref{sec:smoother} we will
discuss the choice of $opt$ in detail.

At the end of this section, let us briefly recall the main idea of multigrid methods. For simplicity, we restrict ourselves to quadratic objective functions $
  f_k(x^{(k)}) = \frac{1}{2} (x^{(k)})^T Q_k x^{(k)} - q_k^T x^{(k)}
$, i.e., equivalently, to systems of linear equations $Q_k x^{(k)} = q_k$.
We start at the finest level $k=j$ and apply some steps of the smoother $opt$ to get an update $x^{(k)}$. We then compute the residuum $r_k = Q_k x^{(k)} - q_k$ for the current iteration. If we now solved the correction equation $Q_{k} e^{(k)} = r_{k}$, we would obtain the exact solution in one step as $x^{(k)}_*=x^{(k)}-e^{(k)}$. This would, however, be as costly as the solution of the original problem. So instead we solve the correction equation on the coarser level, $Q_{k-1} e^{(k-1)} = r_{k-1}$ with $r_{k-1} = I_{k}^{k-1}r_k$, either exactly or approximately by repeating the same procedure on this level. The correction is then interpolated back to level $k$ and a new iteration is computed as $x^{(k)}:=x^{(k)}-I_{k-1}^ke^{(k-1)}$. Finally, we again apply a few steps of the smoother. This is the so-called Correction Scheme (CS) multigrid algorithm. It splits the approximation error into two components, smooth and coarse (relative to the current level); the smooth component is reduced by the smoothing operator, while the coarse component by the coarse grid correction.

\subsection{Truncation}
We now recall the idea of truncation introduced by \mbox{Kornhuber} \cite{Kornhuber1994}; see also \cite{Graser2009}.
Consider the sequence of problems (\ref{eq:2})  obtained by standard finite element discretization of some infinite dimensional problem.
Denote by $\lambda_i^{(k)}$ the finite element basis function associated
with the $i$-th node on discretization level $k$.

Let $x^{(j)}$ be an approximate solution of (\ref{eq:2}) on the finest discretization level and ${\cal A}_j=({\cal A}^\varphi_j(x^{(j)}) \cup {\cal A}^\psi_j(x^{(j)}))$ the corresponding set of active constraints. We will \emph{truncate} the basis functions $\lambda_i^{(j)}$ by putting them equal to zero at the active nodes $\alpha\in {\cal A}_j$:
$$
  \tilde{\lambda}_i^{(j)}(x^{(j)}_\alpha) = 0 \quad\mbox{if}\quad \alpha\in {\cal A}_j\quad\mbox{for all}\ i\,.
$$
The idea is to perform the coarse grid correction in the next iteration of the multigrid method with the truncated basic functions, instead of the original ones. The basis functions for coarse levels are derived from the truncated basis functions on the finest level. Roughly speaking, we consider the active nodes on the finest level to be fixed by homogeneous Dirichlet boundary condition for the next coarse grid correction. After that, we perform the smoothing with the original basis functions; this may change the set of active nodes, so we update the truncated basis and repeat the procedure. Clearly, once the exact active set is detected, the truncated basis does not change any more and the ``truncated'' problem reduces to an unconstrained problem, just as in the classic active set strategy.

The way to implement the truncation depends on the function $f_j$. For instance, if $f_j$ is a quadratic function $f_j(x^{(j)}) = \frac{1}{2} (x^{(j)})^T Q_j x^{(j)} + q_j^T x^{(j)}$ then the truncation of the ``stiffness'' matrix $Q_j$ amounts to putting all rows and columns with indices from ${\cal A}_j$ equal to zero, and analogously for $q_j$. If $f_j$ can still be expressed as a function of $Q_j$ and $q_j$, such as in the problems in Examples~\ref{ex:2} and~\ref{ex:4}, the truncation is analogous. If $f_j$ is a function defined by means of local stiffness matrices $(A_j)_i$, as in the minimal surface problem in Example~\ref{ex:4}, the truncation is performed for all these matrices (again by putting respective rows and columns equal to zero), and the function is evaluated using these truncated matrices.

Let us stress that the truncation is performed explicitly only on the highest discretization level; on coarser grids it is inherited by means of prolongation/restriction operators. This process is again problem dependent. If $f_k$ is defined by means of the global stiffness matrix $Q_k$, then $Q_{k-1} = I_k^{k-1} Q_k I_{k-1}^k,\ k=1,\ldots,j$. If $f_k$ is defined using local stiffness matrices $(A_j)_i$, then $(A_{k-1})_i = I_k^{k-1} (A_k)_i I_{k-1}^k,\ k=1,\ldots,j,$ for all $i$, etc.

In the following, we denote the truncation operator by ${\rm trun}$ and will use notation such as ${\rm trun\,}Q_j$ and ${\rm trun\,}f_j$.

\subsection{Correction scheme for the truncated multigrid for quadratic problems}
In this section we recall the Correction Scheme (CS) version of the truncated monotone multigrid algorithm
introduced by \mbox{Kornhuber} \cite{Graser2009,Kornhuber1994}; see also \cite{krause2001monotone,mandel1984etude,mandel1984multilevel}. We assume that the objective function
$f$ is quadra\-tic, so that
\begin{equation}\label{eq:quadratic}
  f_k(x^{(k)}) = \frac{1}{2} (x^{(k)})^T Q_k x^{(k)} + q_k^T x^{(k)}
\end{equation}
where $Q_k$ are positive definite matrices of size $n_k$ and
$q_k\in\RR^{n_k}$. 
Accordingly, in this section we use the following notation for the smoother:
$$
  x^{(k)}_{\rm new} = opt(Q_k,q_k,\varphi^{(k)},\psi^{(k)},v_k;x^{(k)},\epsilon,\nu)\,.
$$

We introduce two additional restriction operators for the bound
constraints:
\begin{align}
   R^\varphi_k: &\quad (R_k^\varphi y)_i = \max\{y_j \mid j\in{\cal I}_k\cap \mbox{int\,supp\,} \lambda_i^{(k-1)}\},
   \quad i\in {\cal I}_{k-1}\\
   R^\psi_k: &\quad (R_k^\psi y)_i = \min\{y_j \mid j\in{\cal I}_k\cap \mbox{int\,supp\,} \lambda_i^{(k-1)}\},
   \quad i\in {\cal I}_{k-1}\,.
\end{align}

The following is a reformulation of Algorithm 5.10 from \cite{Graser2009} in our notation.
\begin{algorithm}(truncated CS, V-cycle for quadratic problems)
\begin{description}
\item Set $\epsilon,\epsilon_0$. Initialize $x^{(j)}$.
\item for $i=1:niter$
\begin{description}
\item $x^{(j)} := mgm(j,x^{(j)},q_j,\varphi^{(j)},\psi^{(j)})$
\item test convergence
\end{description}
\item end

\item function
    $x^{(k)}=mgm(k,x^{(k)},r_k,\varphi^{(k)},\psi^{(k)})$
\begin{description}
\item if $k=0$
\begin{description}\item $x^{(k)}:= opt({Q}_k,{q}_k,\varphi^{(k)},\psi^{(k)},0;x^{(k)},\epsilon_0,\nu_0)$\hfill (coarsest grid solution)  \end{description} %
\item else%
\begin{description}\item $x^{(k)}:= opt({Q}_k,{q}_k,\varphi^{(k)},\psi^{(k)},0;x^{(k)},\epsilon,\nu_1)$ \hfill (pre-smoothing)
\item $r_k = q_k - Q_k x^{(k)}$\hfill (residuum)
\item $\hat{\varphi}^{(k)} = \varphi^{(k)} - x^{(k)}$
\item $\hat{\psi}^{(k)} = \psi^{(k)} - x^{(k)}$
\item if $k=j$
\begin{description}
\item ${Q}_k := {\rm trun\,} Q_k$ \hfill (matrix truncation)
\item ${r}_k := {\rm trun\,} r_k$ \hfill (residual truncation)
\item $\hat{\varphi}^{(k)}(i)=-\infty$ if $\varphi^{(k)}(i) = x^{(k)}(i)$ \hfill (bounds truncation)
\item $\hat{\psi}^{(k)}(i)=\infty$ if $\psi^{(k)}(i) = x^{(k)}(i)$ \hfill (bounds truncation)
\end{description}
\item end%
\item ${Q}_{k-1} = I_{k}^{k-1} {Q}_{k} I_{k-1}^{k}$ \hfill (coarse grid matrix definition)
\item $r_{k-1} = I_k^{k-1} r_k$\hfill (residal restriction)
\item $\varphi^{(k-1)} = R_k^\varphi (\hat{\varphi}^{(k)})$\hfill (coarse grid bounds)
\item $\psi^{(k-1)} = R_k^\psi (\hat{\psi}^{(k)})$\hfill (coarse grid bounds)
\item $v^{(k-1)} =
    mgm(k-1,0_{n_{k-1}},r_{k-1},\varphi^{(k-1)},\psi^{(k-1)})$\hfill (coarse grid corr.)
\item $x^{(k)} := x^{(k)} + I_{k-1}^k v^{(k-1)}$\hfill (solution update)
\item $x^{(k)}:= opt({Q}_k,{q}_k,\varphi^{(k)},\psi^{(k)},0;x^{(k)},\epsilon,\nu_2)$  \hfill (post-smoothing) \end{description}%
\item end
\end{description}
\end{description}
\end{algorithm}

It is shown in \cite{Graser2009} that the above algorithm converges for any initial iterate.
\subsection{Full approximation scheme truncated multigrid for general problems}
In the above CS algorithm, the coarse grid correction is used to correct the error in the current solution. That is, the solution on the coarse level is just a correction of the fine-level error, not an approximate solution of the original problem. For nonlinear problem, it is useful to solve on every level an approximation of the original problem. This is readily obtained by replacing the restricted residuum $I_k^{k-1}(q_k - Q_k x^{(k)})$ in Algorithm~1 by the ``true'' coarse-level right-hand side corrected by the approximation error, i.e., by $q_{k-1} + (Q_{k-1}I_k^{k-1}x^{(k)} - I_k^{k-1}Q_k x^{(k)})$. The solution update is then changed accordingly. This gives rise to the so-called Full Approximation Scheme (FAS) algorithm.
The FAS version of Algorithm~1 for general nonlinear problems of type (\ref{eq:2}) is given below.

\begin{algorithm}(truncated FAS, V-cycle for nonlinear problems)
\begin{description}
\item Set $\epsilon,\epsilon_0$. Initialize $x^{(j)}$.
\item for $i=1:niter$
\begin{description}
\item $v_j=0_{n_j\times 1}$
\item $x^{(j)} := mgm(f_j,x^{(j)},v_j,\varphi^{(j)},\psi^{(j)})$
\item test convergence
\end{description}
\item end

\item function
    $x^{(k)}=mgm(f_k,x^{(k)},v_k,\varphi^{(k)},\psi^{(k)})$
\begin{description}
\item if $k=0$
\begin{description}\item $x^{(k)}:= opt(f_k,\varphi^{(k)},\psi^{(k)},v_k;x^{(k)},\epsilon_0,\nu_0)$\hfill (coarsest grid solution)  \end{description} %
\item else%
\begin{description}\item $x^{(k)}:= opt(f_k,\varphi^{(k)},\psi^{(k)},v_k;x^{(k)},\epsilon,\nu1)$ \hfill (pre-smoothing)
\item $x^{(k-1)}=I_k^{k-1}x^{(k)}$ \hfill (solution restriction)
\item $g_k=\nabla f_k(x^{(k)})$
\item $\hat{\varphi}^{(k)} = \varphi^{(k)} - x^{(k)}$
\item $\hat{\psi}^{(k)} = \psi^{(k)} - x^{(k)}$
\item if $k=j$
\begin{description}
\item $f_k := {\rm trun\,} f_k$ \hfill (function truncation)
\item $g_k := {\rm trun\,} g_k$ \hfill (gradient truncation)
\item $\hat{\varphi}^{(k)}(i)=-\infty$ if $\varphi^{(k)}(i) = x^{(k)}(i)$ \hfill (bounds truncation)
\item $\hat{\psi}^{(k)}(i)=\infty$ if $\psi^{(k)}(i) = x^{(k)}(i)$ \hfill (bounds truncation)
\end{description}
\item end%
\item $v_{k-1} = I_k^{k-1}v_k + (\nabla f_{k-1}(x^{(k-1)}) - I_k^{k-1}g_k)$ \hfill (correction r.h.s.)
\item $f_{k-1} = I_{k}^{k-1} f_{k} I_{k-1}^{k}$ \hfill (symbolic coarse grid function definition)
\item $\varphi^{(k-1)} = {R}_k^\varphi (\hat{\varphi}^{(k)}) + x^{(k-1)}$\hfill (coarse grid bounds)
\item $\psi^{(k-1)} = {R}_k^\psi (\hat{\psi}^{(k)})+x^{(k-1)}$\hfill (coarse grid bounds)
\item $v^{(k-1)} =
    mgm(f_{k-1},x^{(k-1)},v_{k-1},\varphi^{(k-1)},\psi^{(k-1)})$\hfill (coarse grid corr.)
\item $x^{(k)} := x^{(k)} + I_{k-1}^k (v^{(k-1)}-x^{(k-1)})$\hfill (solution update)
\item $x^{(k)}:= opt({Q}_k,{q}_k,\varphi^{(k)},\psi^{(k)},v_k;x^{(k)},\epsilon,\nu_2)$  \hfill (post-smoothing) \end{description}%
\item end
\end{description}
\end{description}
\end{algorithm}
Note that the coarse grid function definition $f_{k-1} = I_{k}^{k-1} f_{k}$ is only symbolic. As explained in the previous section, it must be defined specifically for each function $f$. Also, the restriction of the solution using operator $I_k^{k-1}$ may not always lead to the best results and one may prefer, for instance, $L_2$ projection; see, e.g., \cite{gross2009convergence}.

It has been observed in \cite[Ex.\,7.3.1]{Graser2009} that the truncated multigrid algorithm may, in the first iterations, be slower than other algorithms (see also Section~\ref{sec:truncate}). This difference may be significant if we are only interested in a low-accuracy solution of the problem. Therefore, in the next section we propose a ``non-truncated'' FAS algorithm for nonlinear problems of type (\ref{eq:2}). This algorithm is a minor generalization of the method proposed in the pioneering paper by Hackbusch and Mittelmann \cite{hackbusch1983}.

\subsection{Full approximation scheme multigrid without truncation}
We now present a Full Approximation Scheme (FAS) version of
Algorithm~2 without truncation of the finite element basis. In order to guarantee convergence of the algorithm, we need to modify the definition of coarse grid constraints. As mentioned above, a similar algorithm has been introduced by Hackbusch and Mittelmann \cite{hackbusch1983}. The difference is in the treatment of the constraints on the coarse levels; while Hackbusch and Mittelmann used active sets, we use the restriction operators defined below.

Let us modify the restriction operators $R_k^\varphi,
R_k^\psi$ introduced in the previous section. The motivation for this
is two-fold: the exact solution of the initial problem should be a
fixed point of the algorithm and a feasible point should remain
feasible after the correction step. For given $x^{(k)}$,  $\varphi$, $\psi$ and some $y\in\RR^{n_k}$ and $i\in {\cal
I}_{k-1}$ the operators are defined as follows
\begin{align}
   \quad (\widetilde{R}_k^\varphi y)_i
   &= \left<
   \begin{aligned}& 0\ \mbox{\it if}\ \max\{(\varphi_j-x^{(k)}_j)\mid
   j\in{\cal I}_k\cap \mbox{int\,supp\,} \lambda_i^{(k-1)}\}=0
\label{eq:5}\\
   &\max\{y_j \mid j\in{\cal I}_k\cap \mbox{int\,supp\,} \lambda_i^{(k-1)}\}\ \mbox{\it otherwise}
   \end{aligned}\right.\\
   \quad (\widetilde{R}_k^\psi y)_i
   &= \left<
   \begin{aligned}& 0\ \mbox{\it if}\ \min\{(\psi_j-x^{(k)}_j)\mid
   j\in{\cal I}_k\cap \mbox{int\,supp\,} \lambda_i^{(k-1)}\}=0\\
   &\min\{y_j \mid j\in{\cal I}_k\cap \mbox{int\,supp\,} \lambda_i^{(k-1)}\}\ \mbox{\it otherwise}\,.
   \end{aligned}\right.\label{eq:6}
\end{align}

\begin{figure}[h]
\begin{center}
 \resizebox{0.5\hsize}{!}
   {\includegraphics{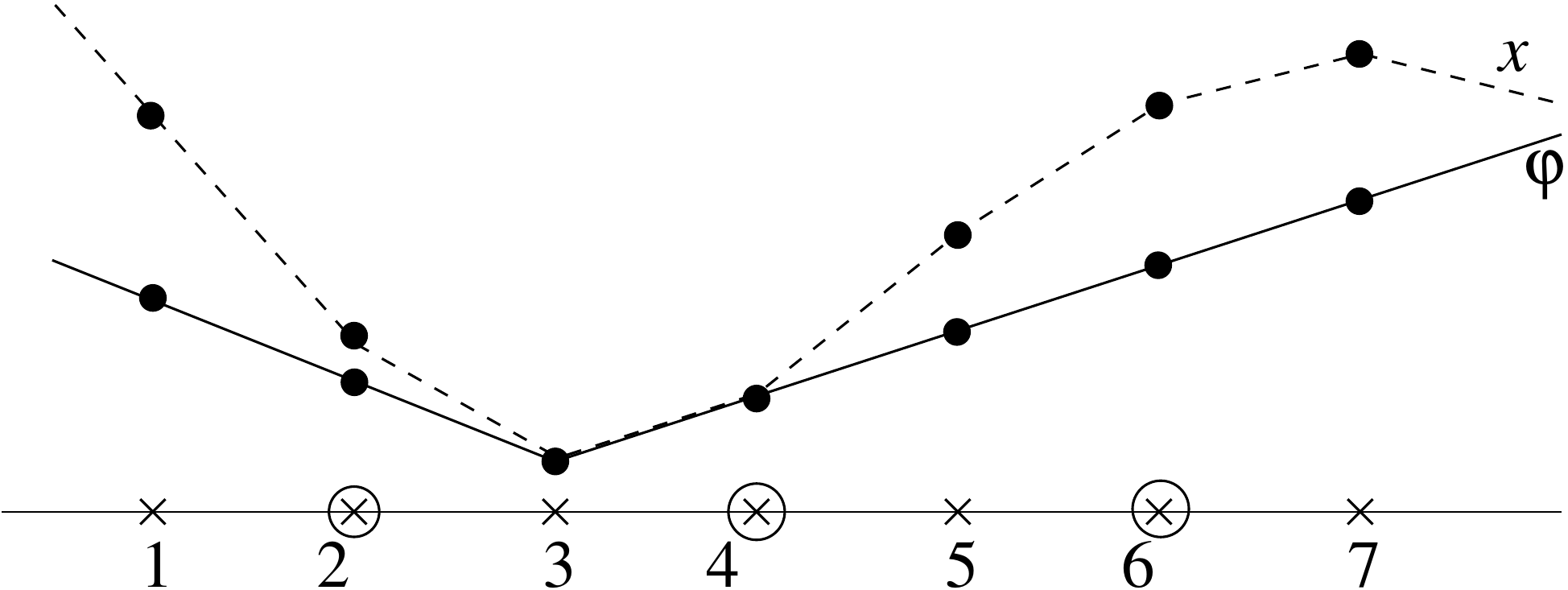}}
  \end{center}
  \caption{A segment of an obstacle and an approximate solution.}\label{fig:con}
\end{figure}

Figure~\ref{fig:con} illustrates how operator $(\widetilde{R}_k^\varphi y)_i$ works.
It depicts a segment of a one-dimensional mesh with three coarse nodes (circles indexed $2,4,6$) and seven fine nodes (crosses indexed $1,\ldots,7$). The constraints are active at nodes $3,4$. For the coarse nodes $2$ and $4$, the first condition applies, as at least one of its neighbours is active. Hence  $(\widetilde{R}_k^\varphi y)_2 =
0$, $(\widetilde{R}_k^\varphi y)_4 =
0$. For node $6$ the second condition in the definition of $(\widetilde{R}_k^\varphi y)_i$ applies, so
$(\widetilde{R}_k^\varphi y)_6 = \max\{y_5,y_6,y_7\}$.  Notice that the operators will not be applied directly to functions $\varphi,\psi$, rather to their modifications; see the details of the FAS algorithm below.

This strategy of handling the coarse-level constraints is a combination of that of Kornhuber (as in Algorithm~1) and the active-set strategy by Hackbush and Mittelmann \cite{hackbusch1983}. We are, however, slightly less conservative than \cite{hackbusch1983}. In their algorithm, both nodes $2$ and $4$ would be considered active in the coarse-level problem and the corresponding value of $x$ would not be allowed to change, unlike in the FAS algorithm below.

\begin{algorithm}(FAS, V-cycle for nonlinear problems)
\begin{description}
\item Set $\epsilon,\epsilon_0$. Initialize $x^{(j)}$.
\item for $i=1:niter$
\begin{description}
\item $q_j = 0_{n_j\times 1}$
\item $x^{(j)} := mgm(j,x^{(j)},q_j,\varphi^{(j)},\psi^{(j)})$
\item test convergence
\end{description}
\item end

\item function
    $x^{(k)}=mgm(k,x^{(k)},q_k,\varphi^{(k)},\psi^{(k)})$
\begin{description}
\item if $k=0$
\begin{description}\item $x^{(k)}:= opt(f_k,\varphi^{(k)},\psi^{(k)},q_k;x^{(k)},\epsilon_0,\nu_0)$  \end{description} %
\item else%
\begin{description}\item $x^{(k)}:= opt(f_k,\varphi^{(k)},\psi^{(k)},q_k;x^{(k)},\epsilon,\nu_1)$ \hfill (pre-smoothing)
\item $x^{(k-1)} = {I}_k^{k-1} x^{(k)}$
\item $q_{k-1} = I_k^{k-1}(q_k-\nabla f_k(x^{(k)})) + \nabla f_{k-1}(x^{(k-1)})$
\item $\hat{\varphi}^{(k)} = \varphi^{(k)} - x^{(k)}$
\item $\hat{\psi}^{(k)} = \psi^{(k)} - x^{(k)}$
\item if $k=j$
\begin{description}
\item $\hat{\varphi}^{(k)}(i)=-\infty$ if $\varphi^{(k)}(i) = x^{(k)}(i)$ 
\item $\hat{\psi}^{(k)}(i)=\infty$ if $\psi^{(k)}(i) = x^{(k)}(i)$
\end{description}
\item end
\item $\varphi^{(k-1)} = \widetilde{R}_k^\varphi
    (\hat{\varphi}^{(k)}) + x^{(k-1)}$\hfill (coarse grid bounds)
\item $\psi^{(k-1)} = \widetilde{R}_k^\psi (\hat{\psi}^{(k)}) +
    x^{(k-1)}$\hfill (coarse grid bounds)
\item $v^{(k-1)} =
    mgm(k-1,x^{(k-1)},q_{k-1},\varphi^{(k-1)},\psi^{(k-1)})$
\item $x^{(k)} := x^{(k)} + I_{k-1}^k (v^{(k-1)} -
    x^{(k-1)})$
\item $x^{(k)}:= opt(f_k,\varphi^{(k)},\psi^{(k)},q_k;x^{(k)},\epsilon,\nu_2)$ \hfill (post-smoothing) \end{description}%
\item end
\end{description}
\end{description}
\end{algorithm}

\section{Equality constraints}\label{sec:eq}
As explained in the introduction, treating general (equality or inequality) constraints by multigrid may be difficult, if not impossible, as we may not be able to find the corresponding restriction operators. This, however, becomes easy in case of a single equality constraint involving all variables, such as
\begin{equation}\label{eq:eq}
  \sum_{i=1}^n x_i = \gamma\,.
\end{equation}
Clearly, this constraint will be present in all discretization levels. We only have to guarantee that, having a feasible point with respect to (\ref{eq:eq}) before the coarse-grid correction step, it will stay feasible after the correction. Let $\gamma_j = \gamma$ and assume that we are on the discretization level $k>0$ and that
\begin{equation}\label{eq:eq1}
  \sum_{i=1}^{n_k} x^{(k)}_i = \gamma_k\,,
\end{equation}
where $n_k$ is the number of variables on the $k$-th level.
Then we define the coarser right-hand side for (\ref{eq:eq}) as
$$
  \gamma_{k-1} = \sum_{i=1}^{n_{k-1}} ({I}_k^{k-1} x^{(k)})_i\,.
$$
If $v^{(k-1)}$ is a solution of the $(k-1)$-level problem, then after the correction step
$x^{(k)} := x^{(k)} + I_{k-1}^k (v^{(k-1)} - {I}_k^{k-1} x^{(k)})$, we obviously get again the equality (\ref{eq:eq1}).

Notice, however, that when we want to combine the equality constraint with the box constraints, we will have to use Algorithm~3 with the \emph{untruncated} restriction operators $(\widetilde{R}_k^\varphi y)_i$ and $(\widetilde{R}_k^\psi y)_i$. The truncated multigrid in Algorithm~2 is not compatible with the equality constraint handled as above.

\section{Smoothing by the steepest descent method}\label{sec:smoother}
The choice of the smoothing method is vitally important for any multigrid algorithm. As our main aim is to avoid second derivatives of the function to be minimized, we have to resort to a first-order optimization method. Moreover, our choice of constrained convex minimization further narrows the choice of available algorithms. We have opted for the gradient projection method. Let us try to justify this choice in the next paragraphs.

\subsection{Steepest descent smoother for unconstrained quadratic problems}
Let us start with an unconstrained convex quadratic problem
\begin{equation}\label{eq:sm1}
  \min_x \frac{1}{2} x^T Q x - q^Tx
\end{equation}
or, in other words, with a linear system
\begin{equation}\label{eq:sm2}
   Q x = q
\end{equation}
and the classic V-cycle multigrid algorithm. If things wouldn't work here, we can hardly expect them to work in the more general setting. One of the most popular smoothers in this case is the Gauss-Seidel (GS) iterative method. Rightly so, its very definition shows that it solves the  equations locally, one by one, performing thus the local smoothing of the approximate solution of (\ref{eq:sm2}).

Looking at the optimization formulation (\ref{eq:sm1}) of the problem, we can also consider the ``most basic'' optimization algorithm, the steepest descent (SD) method with line search. Can this be a good smoother? This question was analyzed, e.g., by McCormick \cite{mccormick1985multigrid} who showed that the steepest descent method with exact or ``slightly inexact'' line search has indeed smoothing properties, as required for the convergence of the standard V-cycle. In case of exact line search, McCormick also gives an explicit bound on the convergence of the V-cycle using the steepest descent method as a smoother. This bound is, however, as many such theoretical bounds, overly pessimistic and far away from the real behaviour of the method (giving estimates of convergence speed such as 0.9995).

We have therefore performed a small experiment with the goal of testing the smoothing property of the steepest descent method and comparing it to the Gauss-Seidel method. Let us introduce some notation. Denote the exact solution on level $h$ by $(x^*)^h$. Consider two discretizations of the underlying differential equation, one on the fine level parameterized by $h$ and one on the coarse level $2h$. Let $P$ be a prolongation operator from the coarse to the fine level. It is well known (see, e.g., \cite{mccormick1985multigrid}) that the $Q$-orthogonal projector on the range of $P$ can be written as
$$
  S^h = P(P^T Q^h P)^{-1}P^T Q
$$
while
$$
  T^h = I-S^h
$$
is the projector onto the $Q$-orthogonal complement of the range of $P$. With these two projectors, the energy norm of the error  $e^h=(x^*)^h-x^h$ of an approximate solution $x^h$ to (\ref{eq:sm2}) satisfies
$$
  \|e^h\|^2_{Q^h} = \|S^h e^h\|^2_{Q^h} + \|T^h e^h\|^2_{Q^h}
$$
whereas $S^h e^h$ and $T^h e^h$ are the ``low-frequency" and the ``high-frequency" components of the error. It is the goal of the smoother to reduce $T^h e^h$ quickly in a few initial (possibly just 1--2) iterations.

In our experiment, the underlying problem was the Poisson problem on a unit square discretized by standard quadrilateral bilinear finite elements. We consider a 32 by 32 fine grid and 16 by 16 coarse grid; the condition number of the fine-grid matrix was 400. The restriction operator is the standard full weighting operator, see, e.g., \cite{briggs2000multigrid}. We consider the steepest descent method
$ x^h_{k+1} = x^h_k - s_k (Q^h x^h_k - q^h)$
with exact step length $s_k$ and a method with inexact line search, the details of which are given in the next section.
We start the process with a randomly generated vector $x_0^h$ using the MATLAB function {\tt randn}. Assuming that the right-hand side, and thus the solution, is smooth, we will get significant components in both $S^h e^h_0$ and $T^h e^h_0$.

\begin{figure}[h]
\begin{center}
 \resizebox{0.45\hsize}{!}
   {\includegraphics{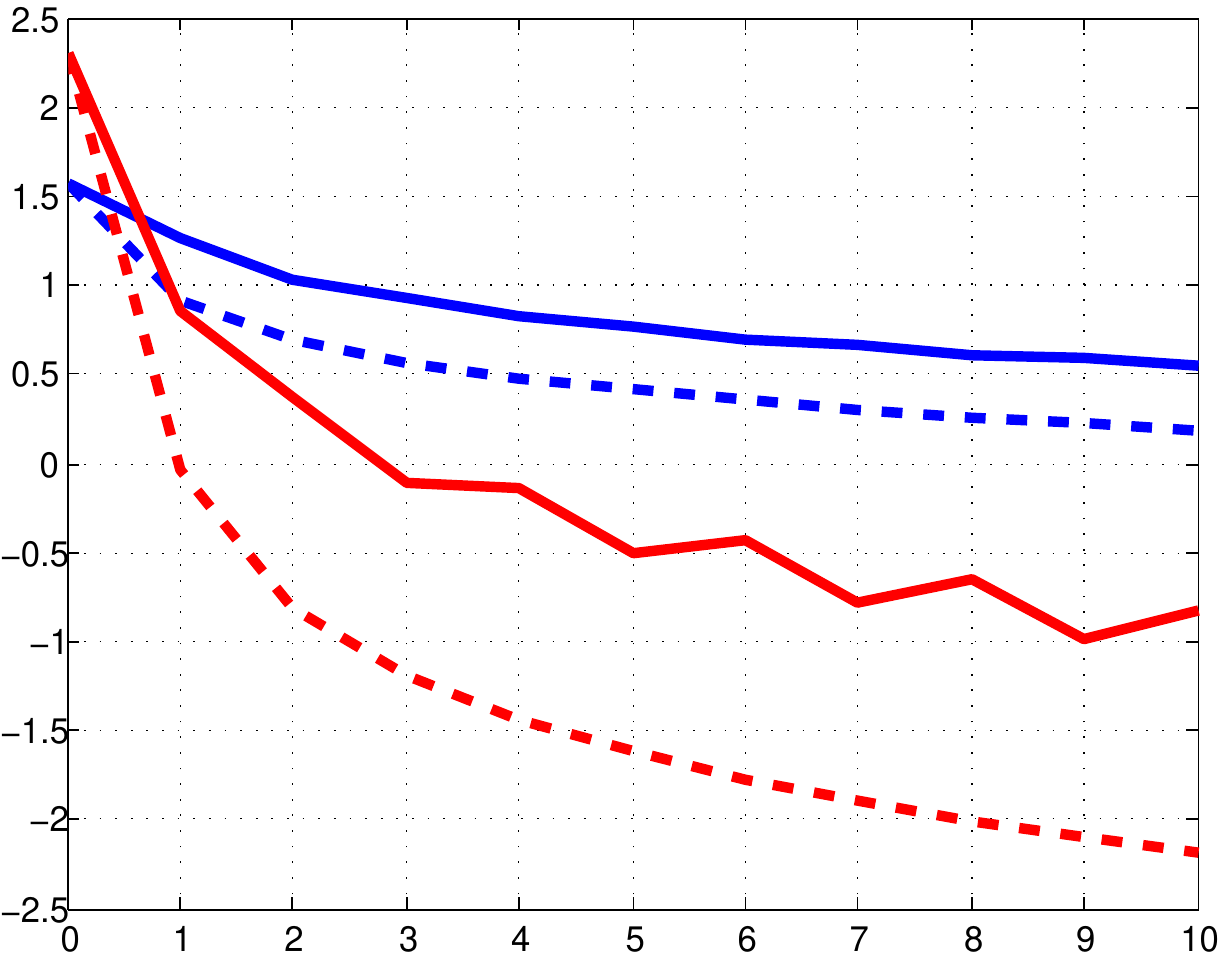}}\quad
 \resizebox{0.45\hsize}{!}
   {\includegraphics{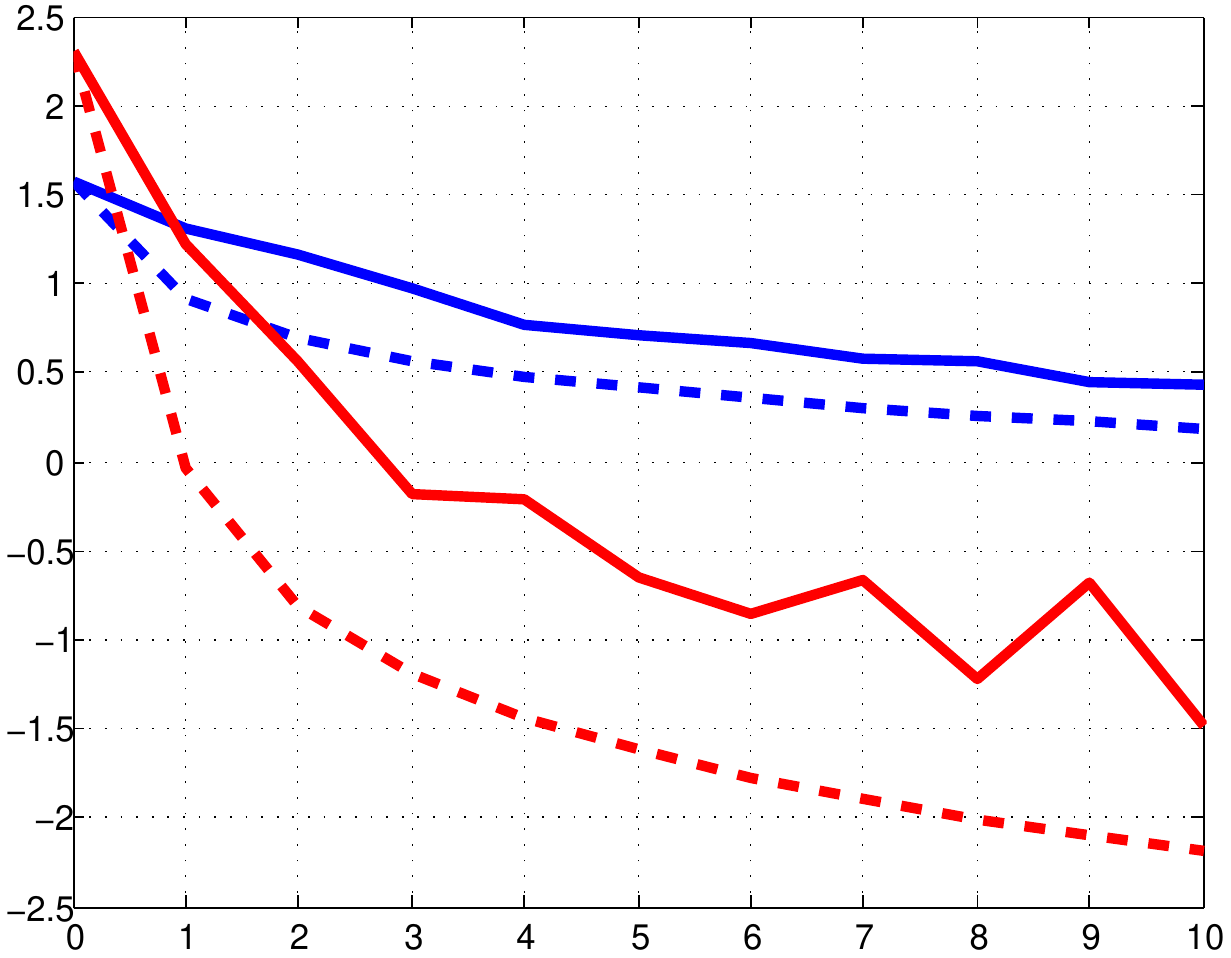}}
  \end{center}
  \caption{First 10 iterations of SD (full line) and GS (dashed line) methods. Blue line depicts low frequency error $\|S^h e^h\|^2_{A^h}$, red line the high frequency error $\|T^h e^h\|^2_{A^h}$. SD with exact line search is on the left, with inexact line search on the right.}\label{fig:1}
\end{figure}

\begin{figure}[h]
\begin{center}
 \resizebox{0.45\hsize}{!}
   {\includegraphics{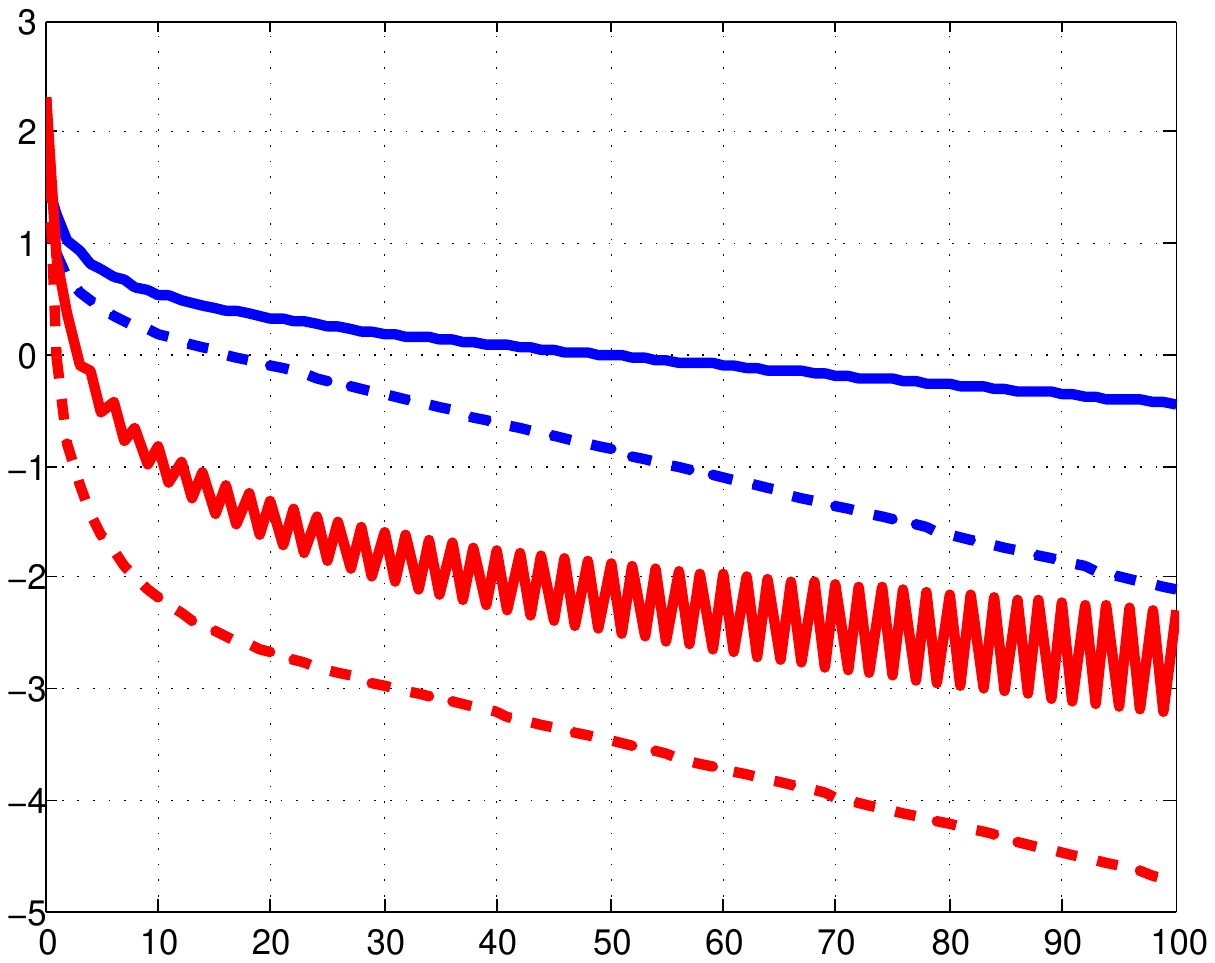}}\quad
 \resizebox{0.45\hsize}{!}
   {\includegraphics{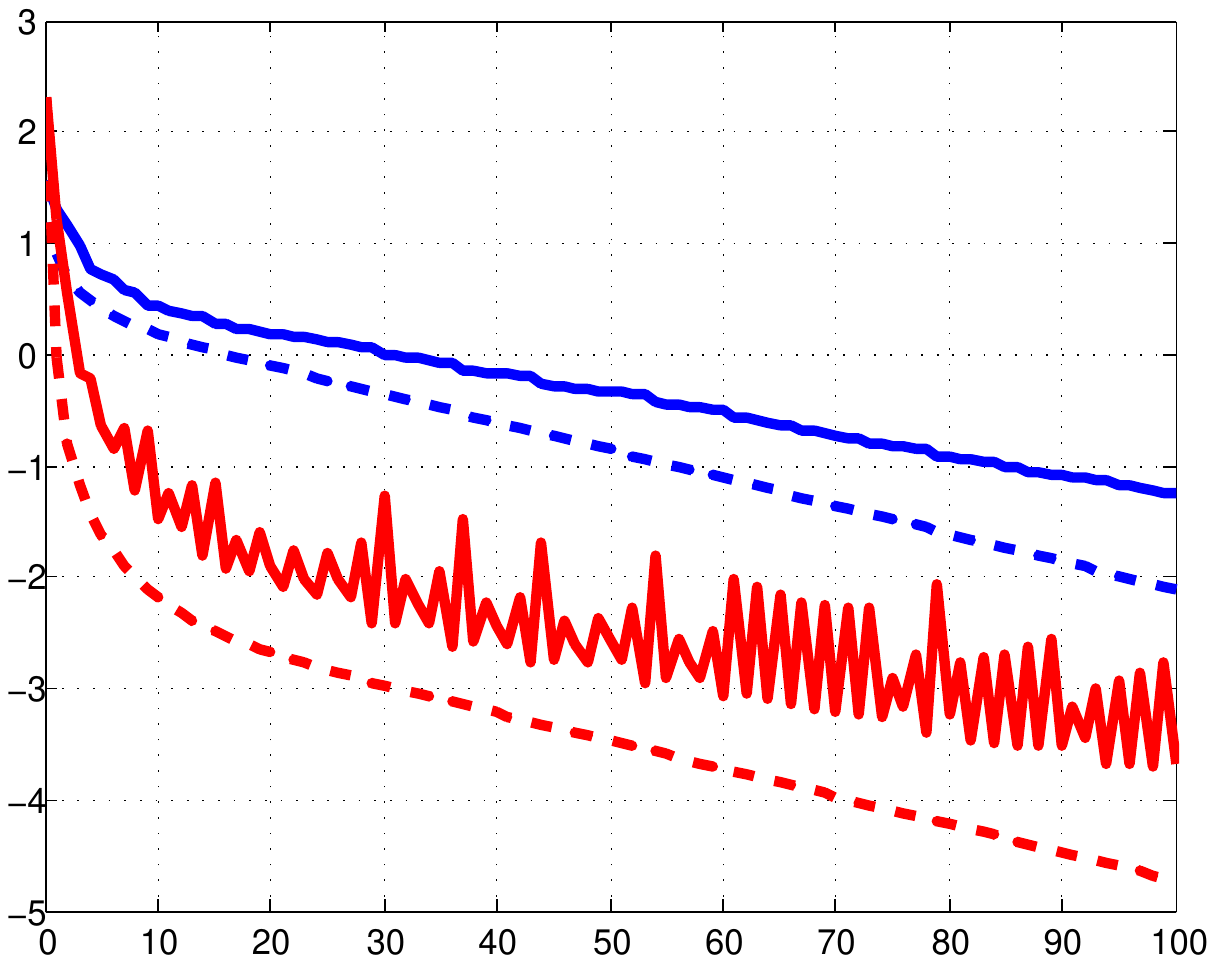}}\\
 \resizebox{0.45\hsize}{!}
   {\includegraphics{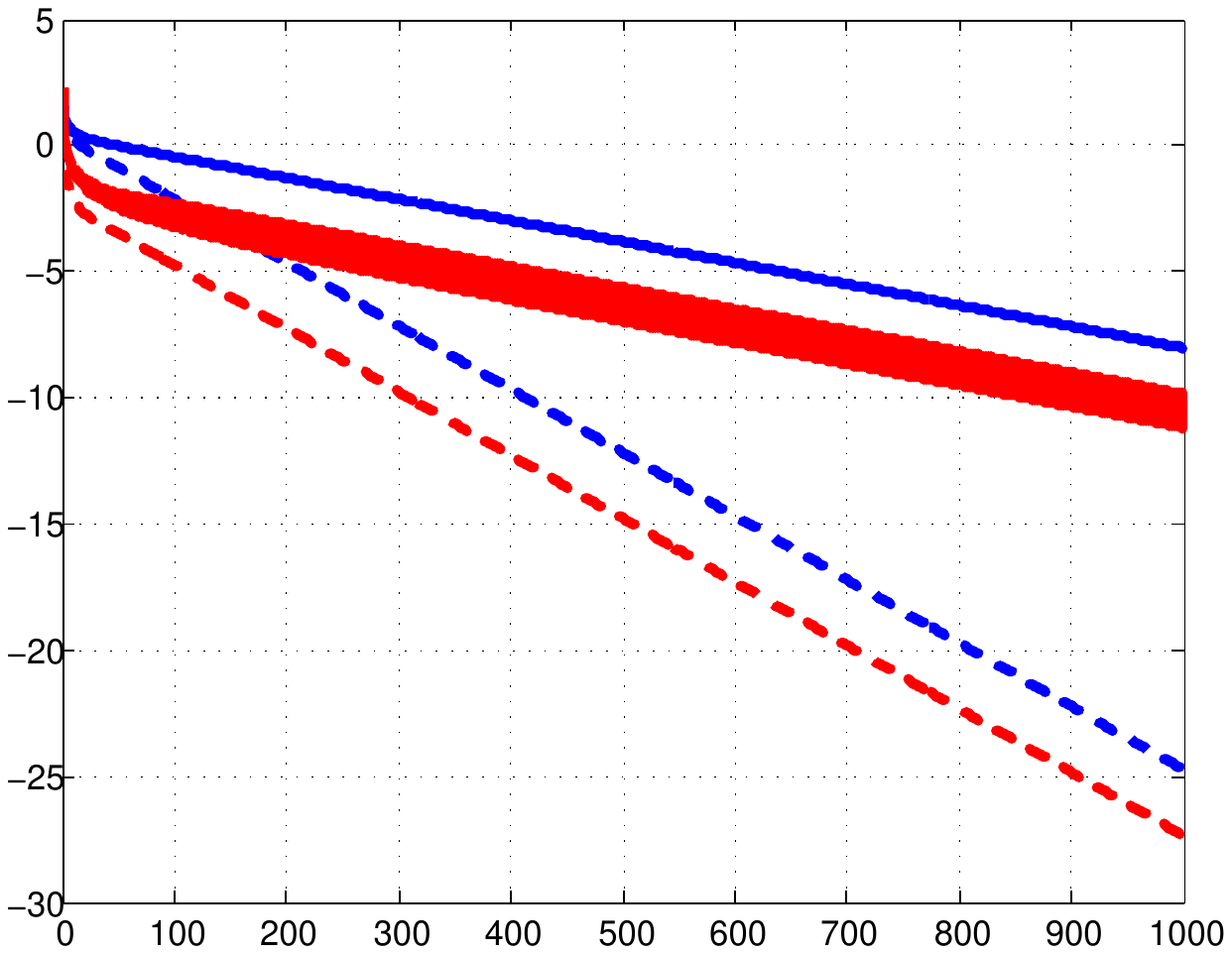}}\quad
 \resizebox{0.45\hsize}{!}
   {\includegraphics{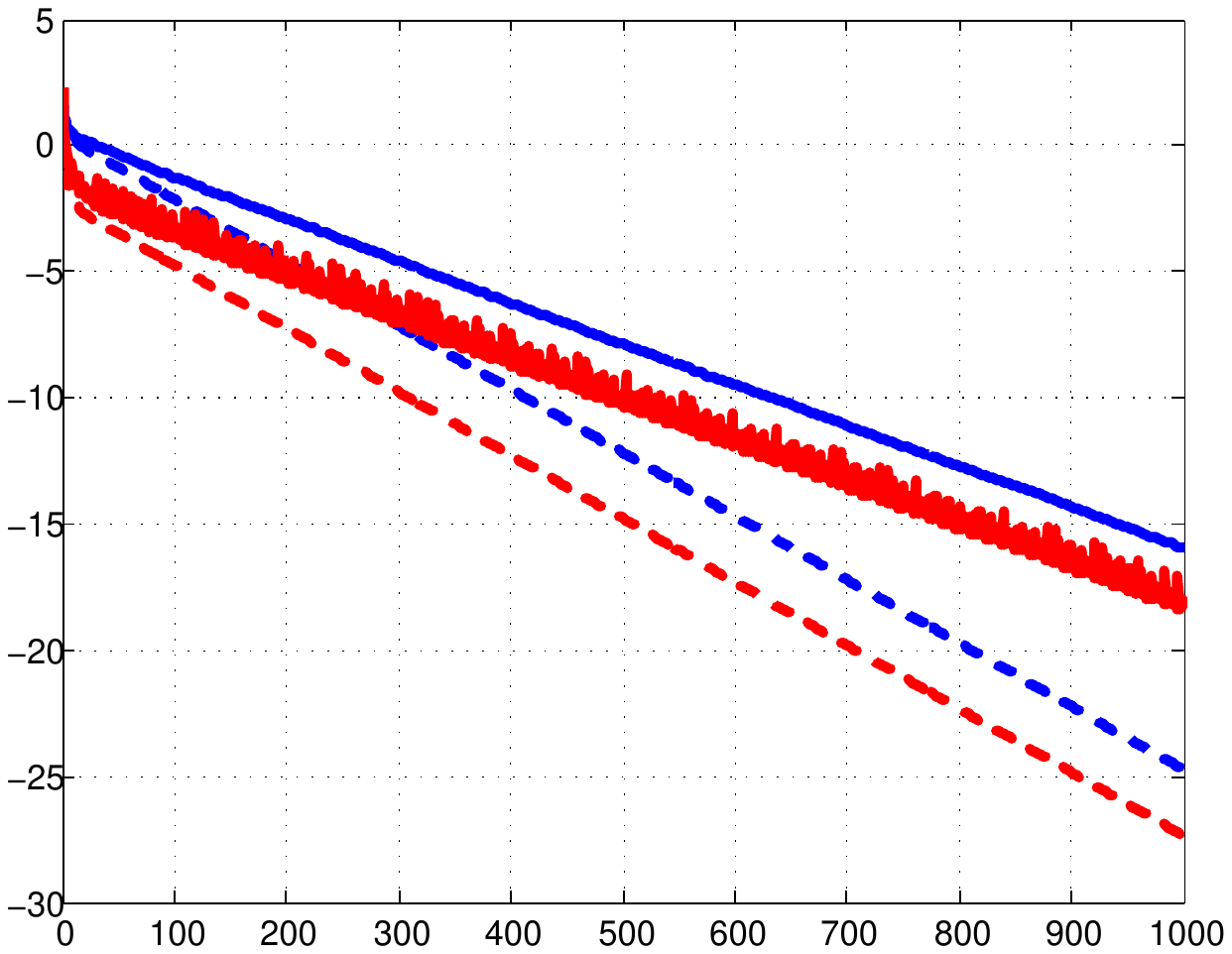}}
  \end{center}
  \caption{First 100 (top) and 1000 (bottom) iterations of SD (full line) and GS (dashed line) methods. Blue line depicts low frequency error $\|S^h e^h\|^2_{A^h}$, red line the high frequency error $\|T^h e^h\|^2_{A^h}$. SD with exact line search is on the left, with inexact line search on the right.}\label{fig:2}
\end{figure}

We first performed 10 iterations of SD and monitored the energy norm of both components of the error, $\|S^h e^h\|^2_{Q^h}$ and  $\|T^h e^h\|^2_{Q^h}$. In Figure~\ref{fig:1}, these are depicted by the full line, the low frequency $\|S^h e^h\|^2_{Q^h}$ in blue and the high-frequency $\|T^h e^h\|^2_{Q^h}$ in red. Both values are given in logarithmic scale. The left-hand figure is for SD with exact line search while the right-hand side one is for the inexact line search. Also in Figure~\ref{fig:1}, we plot these values for the Gauss-Seidel method; these are denoted by the dashed lines.
We can clearly see the smoothing effect of both methods in the first iterations when the red lines quickly drop by orders of magnitude. Although the steepest descent method is not as efficient a smoother as Gauss-Seidel, it still does a good job.
After the initial iterations, the smoothing effect slows down and both errors descent proportionally. This can be better seen in Figure~\ref{fig:2} where we show the error after 100 and after 1000 iterations of both methods.
In this figure, we can clearly see the typical zig-zagging of the SD method with exact line search present in the high frequency error. We can also see that inexact line search with its random element breaks this regular zig-zagging and, in effect, makes the method significantly faster.

We could certainly consider other first-order method, for instance the nonlinear conjugate gradients. However, as we will see in the next section and later in the numerical experiments, only a very few (1--5) iterations of the steepest descent method suffice to guarantee a good behaviour of the multigrid algorithm, and there is thus no need for any more sophisticated first-order algorithms.

\subsection{Line search}
We have seen in the previous section that exact line search does not bring any significant benefit to the steepest descent method. Moreover, our ultimate goal is to use the projected gradient version of the method for bound-constrained nonlinear convex problems.

A popular---and efficient---choice of the step length is the Barzilai-Borwein method \cite{dai2005projected}. This, however, leads to a possibly non-monotonic progress of the error. Since we would like to use a very small \emph{fixed} number of SD steps, this method is not suitable. Moreover, its projected gradient version is not fully understood and may lead to a standard Armijo step \cite{dai2005projected}.

Because our problem is convex, the line search can be based solely on the gradient information. For an unconstrained problem
$$
  \min_x f(x)
$$
with $f$ smooth and convex, we have opted for the following simplified version of Wolfe's method.

\begin{algorithm}(Steepest descent method with gradient-based line search)
\\
Given an approximate solution $x$, do until convergence:
\begin{description}
\item $x_{\rm new}= x-s\nabla f(x)$ with $s$ computed by Algorithm 5.
\end{description}
\end{algorithm}
\begin{algorithm}(Gradient-based line search for unconstrained problems)
\\
Given an approximate solution $x$. Choose $s>0$ and $c>1$.
\begin{description}
\item 1.    $\gamma = -(\nabla f(x))^T\, \nabla f(x-s\nabla f(x))$
\item 2. if $\gamma<0$
\begin{description}
\item 2.1 do until $\gamma^+>0$
\begin{description}
\item 2.1.1 $s:= cs$
\item 2.1.2. $\gamma^+ = -(\nabla f(x))^T\, \nabla f(x-s\nabla f(x))$
\end{description}
\item 2.2. $s:= \frac{1}{c}s$
\end{description}
\item[\quad ] else
\begin{description}
\item 2.3. do until $\gamma^+<0$
\begin{description}
\item 2.3.1. $s:= \frac{1}{c}s$
\item 2.3.2. $\gamma^+ = -(\nabla f(x))^T\, \nabla f(x-s\nabla f(x))$
\end{description}
\end{description}
\item[\quad ] end
\item 3. return current $s$

\end{description}
\end{algorithm}

Clearly, $\gamma$ (or $\gamma^+$) is the directional derivative of $f$ in the steepest descent direction $-\nabla f(x)$ at the trial point $x-s\nabla f(x)$.

\begin{lemma}\label{th:LS1}
Algorithm 5 is well-defined and finishes in a finite number of steps. At the new point, we have
$f(x_{\rm new})<f(x)$ and $\gamma^+ < -(\nabla f(x))^T\, \nabla f(x_{\rm new})<0$.
\end{lemma}
\begin{Proof}
The claim follows immediately from the convexity of $f$. The algorithm stops with the value of $s$ for which $\gamma^+<0$ and such that for the step length $cs$, the derivative would change sign to $\gamma^+>0$. The loop in 2.1 stops when we encounter a positive $\gamma^+$, so we have to return one step back in 2.2.
\end{Proof}

Let us now move to a generalization of Algorithm 4 for convex bound-constrained problem (\ref{eq:1}).
For a given feasible $x$, let us recall the definition of the set of active indices by
$$
  {\cal A}(x) = \left\{i\mid x_i=\varphi_i\ \mbox{or}\ x_i=\psi_i\right\}\,.
$$
For this $x$, we further introduce an operator  $[\,\cdot\,]_{{\cal A}(x)} : \RR^n\to\RR^n$ defined component-wise by
$$
  ([\,h\,]_{{\cal A}(x)})_i =\left<\begin{aligned}
  0 &\ \mbox{if}\ i\in{\cal A}(x)\\
  h_i&\ \mbox{otherwise.}\end{aligned}\right.
$$
Finally, let $[\,\cdot\,]_\Omega$ denote the (in this case trivial) projection on the feasible set.

\begin{algorithm}(Gradient projection method with gradient-based line search)\\
Given a feasible approximate solution $x\in\Omega$, do until convergence
\begin{description}
\item $x_{\rm new} = [x-s\nabla f(x)]_\Omega$ with $s$ computed by Algorithm 7.
\end{description}
\end{algorithm}

\begin{algorithm}(Gradient-based line search for constrained problems)\\
Given a feasible approximate solution $x\in\Omega$. Choose $s>0$ and $c>1$.
\begin{description}
\item 1.    $x^+ = [x-s\nabla f(x)]_\Omega$,\quad $\gamma = -(\nabla f(x))^T\, [\nabla f(x^+)]_{{\cal A}(x^+)}$
\item 2. if $\gamma<0$
\begin{description}
\item 2.1 do until $\gamma^+>0$
\begin{description}
\item 2.1.1 $s:= cs$
\item 2.1.2. $x^+ = [x-s\nabla f(x)]_\Omega$
\item 2.1.3. $\gamma^+ = -(\nabla f(x))^T\, [\nabla f(x^+)]_{{\cal A}(x^+)}$
\end{description}
\item 2.2. $s:= \frac{1}{c}s$
\end{description}
\item[\quad ] else
\begin{description}
\item 2.3. do until $\gamma^+<0$
\begin{description}
\item 2.3.1. $s:= \frac{1}{c}s$
\item 2.3.2. $x^+ = [x-s\nabla f(x)]_\Omega$
\item 2.3.3. $\gamma^+ = -(\nabla f(x))^T\, [\nabla f(x^+)]_{{\cal A}(x^+)}$
\end{description}
\end{description}
\item[\quad ] end
\item 3. return current $s$
\end{description}
\end{algorithm}

\begin{figure}[h]
\begin{center}
 \resizebox{0.48\hsize}{!}
   {\includegraphics{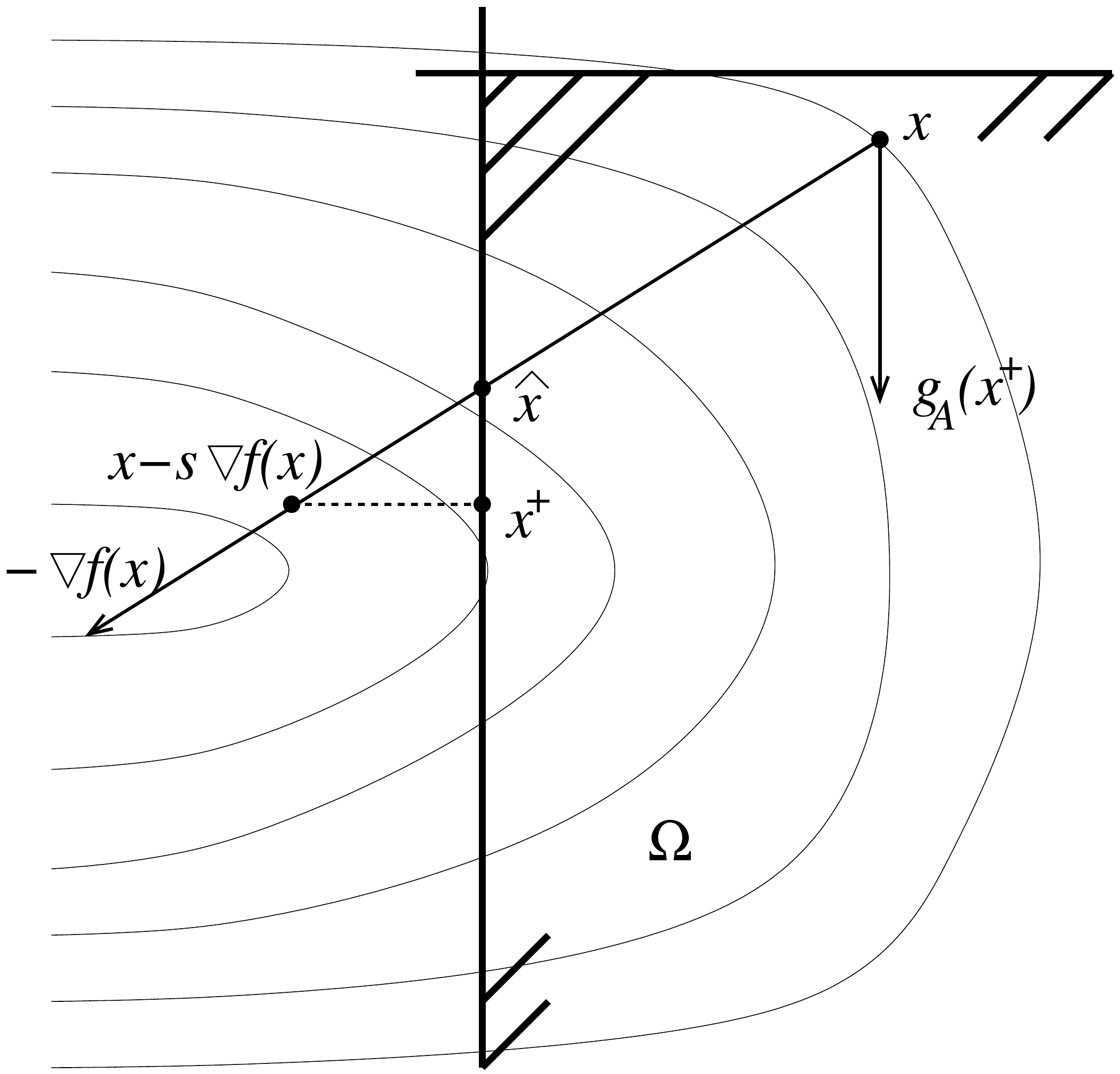}}\quad
 \resizebox{0.48\hsize}{!}
   {\includegraphics{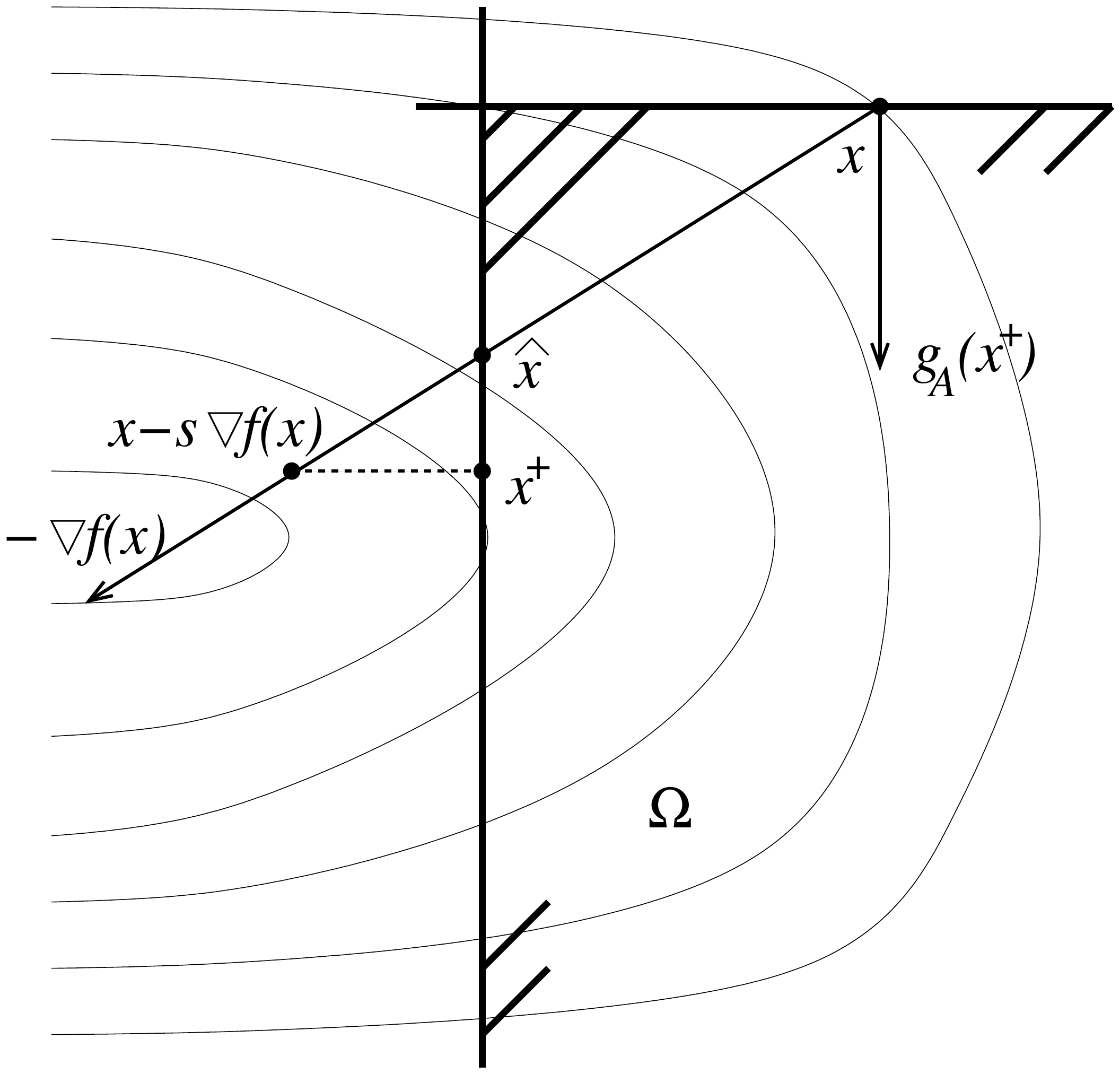}}
  \end{center}
  \caption{Gradient-based line search.}\label{fig:grad}
\end{figure}

\begin{lemma}\label{th:LS2}
Algorithm 7 is well-defined and finishes in a finite number of steps. At the new point, we have
$f(x_{\rm new})<f(x)$ and $\gamma^+ < -(\nabla f(x))^T\, \nabla f(x_{\rm new})<0$.
\end{lemma}
\begin{Proof}
For simplicity, we assume that $\gamma<0$, so we are in the 2.1.--2.2.\ branch of the algorithm. The other case would be handled analogously.
Assume first that no constraints are active at the initial point $x$, i.e., ${\cal A}(x)=\emptyset$. Once we have computed a new point $x^+$ and found the active set ${\cal A}(x^+)$, we can split the search direction $-\nabla f(x)$ into two vectors:
$$
  -\nabla f(x) = g_A(x) + g_N(x)
$$
where
$$
  (g_A(x))_i = \left<\begin{aligned}
0\quad &\mbox{if}\ i\in{\cal A}(x^+)\\
-(\nabla f(x))_i\quad &\mbox{otherwise}
\end{aligned}\right.
$$
and complementary for $g_N(x)$. Both $g_A(x)$ and $g_N(x)$ are still descent directions, in particular
$$
  -(\nabla f(x))^T g_A(x)<0\,.
$$
Let
$$
  \hat{x} := x- \hat{s}\nabla f(x)
$$
with $\hat{s}$ chosen such that ${\cal A}(\hat{x})={\cal A}(x^+)$ (see Figure~\ref{fig:grad}, left). If $-(\nabla f(\hat{x}))^T g_A(x)<0$ then the function is still descending at $\hat{x}$ in direction $g_A(x)$ and we do a line search in this direction, i.e., along the manifold defined by active indices. Due to convexity of $f$, we either have to reach a (finite) point at which the directional derivative changes its sign or when we hit a new constraint and the active set changes; in the latter case we repeat the above argument with the new active set. If $-(\nabla f(\hat{x}))^T g_A(x)\geq 0$ then we know that we went too far, the algorithm stops and returns the previous trial point. The rest follows from Lemma~\ref{th:LS1}.

Now if ${\cal A}(x)\not=\emptyset$ then either ${\cal A}(x^+)={\cal A}(x)$ and the above argument applies or some constraints from ${\cal A}(x)$ are released at $x^+$ (Figure~\ref{fig:grad}, right). But this would mean that the search direction goes away from these constraints, they can be ignored and, again, the above arguments apply.

In the (unlikely) case ${\cal A}(x^+)=\{1,2,\ldots,n\}$ when all components of the trial point $x^+$ are active, we have $\gamma^+=0$ and the algorithm stops and returns the previous trial point.
\end{Proof}

Our preferred choice of the parameters is $c=2$ and $s=1$ in the first iteration. In the following iterations of the multigrid V-cycle, the parameter $s$ is chosen as the final one from the previous call of gradient Algorithms 5 or 7.

\section{Numerical experiments}
In this section we presents results of our numerical experiments using examples collected from the literature. We will start with a quadratic problem, in order to see the influence of non-linearity on the behaviour of the multigrid method.

For each example we will present the computed asymptotic rate of convergence for different numbers of the smoothing steps, namely for 2,4,6,8,10 smoothing steps (half of them in the pre-smoothing phase, half in the post-smoothing). In the same table, we will give the number of function and gradient calls on the finest level. We will also present a comparison (in terms of function/gradient calls) with one of the most efficient codes for these problems, the L-BFGS-B by Morales and Nocedal \cite{morales2011remark}.

All examples are defined on the square $\Omega=(0,1)^2$ in the infinite dimensional setting. We then use regular meshes of square finite elements with bilinear basis functions for their discretization.
The initial coarsest mesh (refinement level 0) consists of four elements. We apply $j=8$ uniform refinement steps to get 9 embedded finite element meshes. That means that, on refinement level $k$, $k=0,1,2,\ldots$, we have $4^{k+1}$ finite elements and $(2^{k+1}-1)^2$ interior nodes, with the finest mesh having 262\,144 finite elements and 261\,121 interior nodes.
The prolongation operators $I_{2h}^h$ are based on the nine-point interpolation scheme defined by the stencil $\begin{pmatrix}
\frac{1}{4}&\frac{1}{2}&\frac{1}{4}\\
\frac{1}{2}&1&\frac{1}{2}\\
\frac{1}{4}&\frac{1}{2}&\frac{1}{4}\\
\end{pmatrix}$. We use the full weighting restriction operators defined by $I_h^{2h}= \frac{1}{4}(I_{2h}^h)^T$; see, e.g., \cite{hackbusch1986multigrid}.
The initial point in all experiments was set to a zero vector, even if this was infeasible. The problems on the coarsest level were solved by the same iterative method used as a smoother, however, with high accuracy. In particular, the parameters $\epsilon_0$ and $\nu_0$ in Algorithms 1--3 were set to $\epsilon_0=10^{-9}$ and $\nu_0=10000$. Recall that $\epsilon_0$ controls the norm of the scaled gradient or the KKT conditions for constrained problems and $\nu_0$ is a bound on the number of iterations; because our coarsest problems have very low dimensions, this iteration bound was never reached.

The approximate asymptotic convergence rate is computed as $\displaystyle\left(\frac{e_k}{e_2}\right)^{\frac{1}{k-1}}$, where $e_k=\|x^*-x_k\|$, $x_k$ is the last iteration before the algorithm stops and $x^*$ is the ``exact'' solution as computed either by L-BFGS-B or by the gradient projection method with high accuracy. Every norm in the examples is the Euclidean norm.

As a measure of efficiency of the algorithms we present the number of function and gradient evaluations on the finest level. Of course the work on the lower levels is not free, even though the fine level is dominant, in particular for the nonlinear problems. For instance, in our implementation of the minimum surface problem (which can certainly be improved), the function and gradient evaluation on a fine level was 16 times more expensive than on the coarser level. In most examples, the number of function/gradient evaluations on a a single coarser level is about the same as on the finest level. In the most favourable case when the computational complexity of one function and gradient evaluation is linear in the number of variables, we would expect the effort to compute the function on the coarser level to be about one-quarter the effort on the finer level. This claim is confirmed by a remark concluding Example~\ref{ex:2} presenting exact timings on all levels.

All algorithms were implemented in MATLAB. The interface to the L-BFGS-B code was provided by Stephen Becker\footnote{http://www.mathworks.co.uk/matlabcentral/fileexchange/35104-lbfgsb--l-bfgs-b--mex-wrapper}. For all experiments we used a laptop with Intel Core i7-3570 CPU M 620 at 2.67GHz with 4GB RAM, and MATLAB version 8.0.0 (2012b) running in 64 bit Windows 7.

\subsection{Example: quadratic obstacle problem}\label{ex:1}
\begin{figure}[h]
\begin{center}
 \resizebox{0.48\hsize}{!}
   {\includegraphics{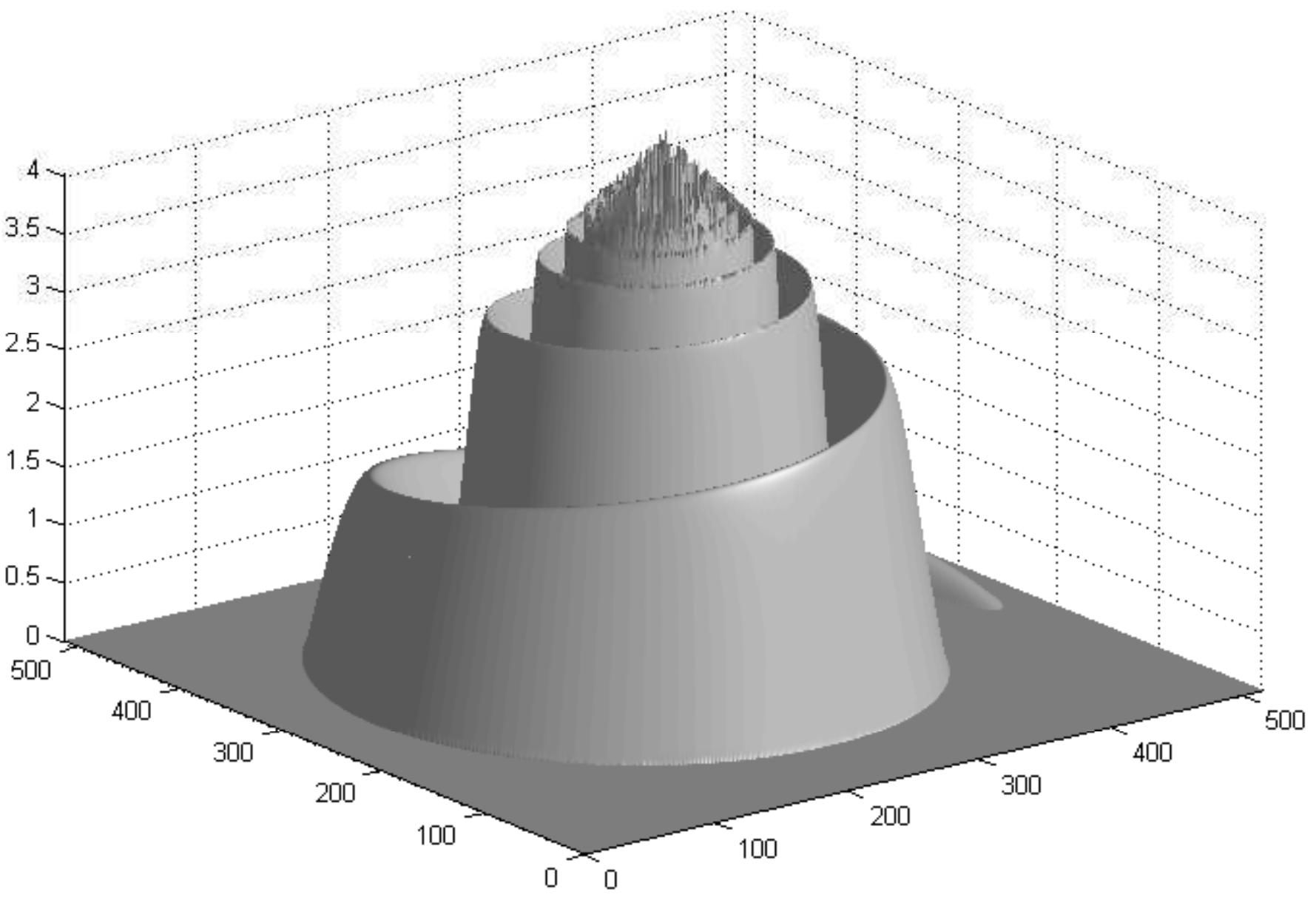}}\quad
 \resizebox{0.48\hsize}{!}
   {\includegraphics{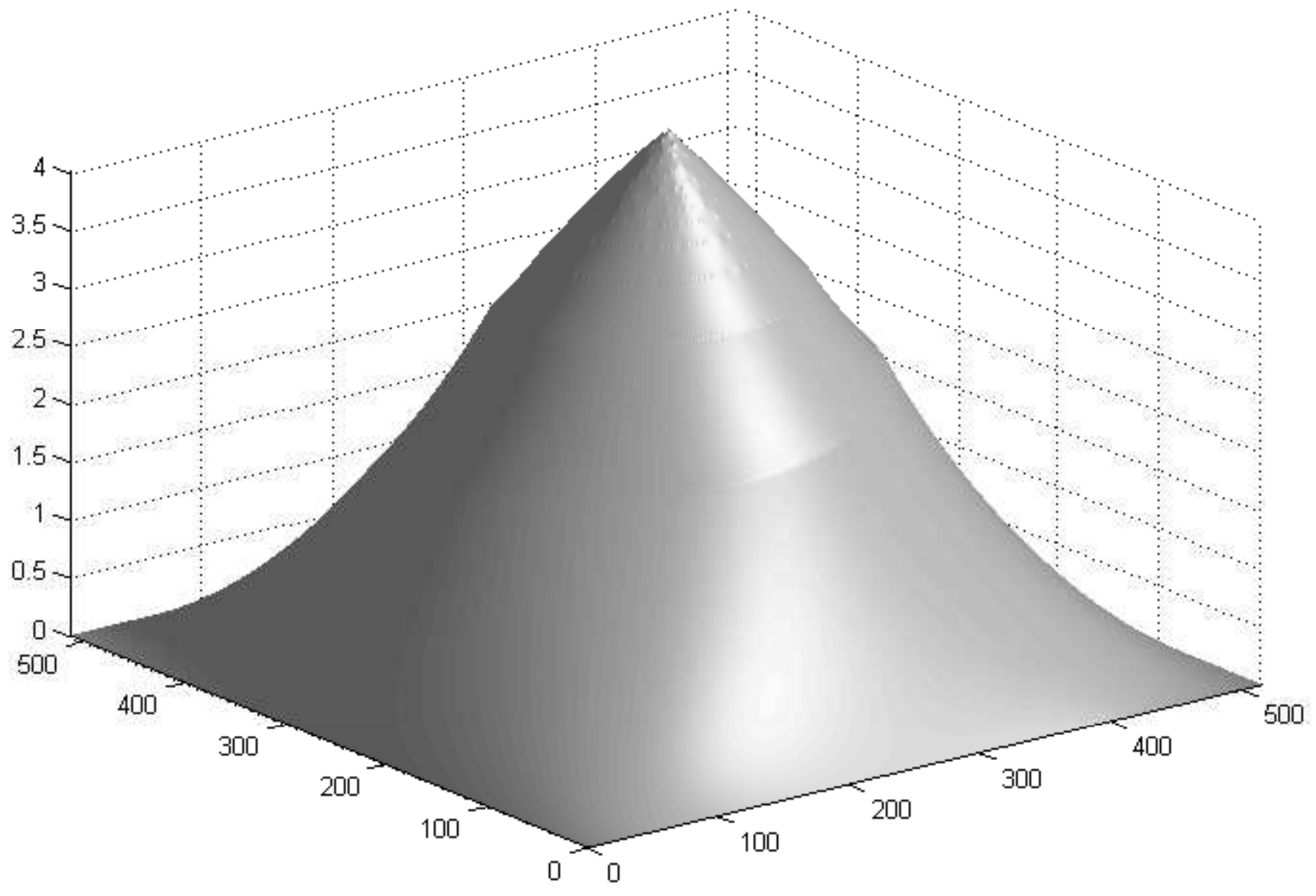}}
  \end{center}
  \caption{Example~\ref{ex:1}, eight refinement levels, obstacle (left) and solution (right).}\label{fig:ex1}
\end{figure}

Let us start with the ``Spiral problem" from \cite{Graser2009}. This is a quadratic optimization problem resulting from the Laplace equation in $\Omega\subset \RR^2$:
\begin{align*}
  &\min_{u\in H^1_0(\Omega)} {\cal J}(u):= \frac{1}{2}\int_\Omega \|\nabla u\|^2\;dx - \int_{\Omega} Fu\;dx\\
  &\mbox{subject to}\nonumber\\
  &\qquad \varphi \leq u \leq \psi,\quad \mbox{a.e.\ in}\ \Omega\,,
\end{align*}
with $F\in L^2(\Omega)$.
We will use the spiral obstacle, as proposed in \cite[\textsection 7.1.1]{Graser2009}:
$$
  \varphi(x(r,\phi)) = \sin(2\pi/r + \pi/2 - \phi)+\frac{r(r+1)}{r-2} - 3r + 3.6,\quad r\not=0\,,
$$
and $\varphi(0)=3.6$ with polar coordinates $x(r,\phi)=re^{i\phi}$. The upper bound function $\psi$ is set to infinity and the right-hand side function $F$ to zero. The obstacle function is illustrated in Figure~\ref{fig:ex1} (left), together with the solution of the problem for 8 refinement levels.

Table~\ref{tab:1} presents the results of the numerical experiments. It shows the asymptotic rate of convergence (an average from the last 3--5 iterations) and the number of evaluations of the objective function and its gradient on the finest level only. The results are presented for 4--8 refinement levels and refer to Algorithm~1 with $\nu$ pre-smoothing and $\nu$ post-smoothing steps (GP-$\nu$). For comparison, we also show the convergence rate for the algorithm with the projected Gauss-Seidel (GSP) smoother (one pre- and one post-smoothing step). The last row of the table presents the numbers of function evaluations when we solved the finest level problem directly by the gradient projection method. The Gauss-Seidel smoother is clearly superior to GP-1 but its convergence rate can be reached by several GP smoothing steps. This increased number of GP steps is size-dependent, as indicated by the examples. However, the results also suggest that the increased number of GP steps is not needed, as the smallest number of function/gradient evaluations is obtained with 1 or 2 smoothing steps. Notice also that one GSP step is much more CPU expensive than one GP step, at least in MATLAB implementation which allows for the vectorization of the GP step. The reason for this is that in GSP the projection must be performed for each variable separately, after its update by the Gauss-Seidel inner iteration, while in the GP the whole vector is projected at once.
\begin{table}[htbp]
  \centering
  \caption{Example~\ref{ex:1}, asymptotic rate of convergence and number of top-level function evaluations for 4--8 refinement levels. Here ``vars" stands for the number of variables, GSP for the (1,1) V-cycle with Gauss-Seidel method with projection; GP-$\nu$ for a $(\nu,\nu)$ V-cycle with the GP smoother; ``GP only'' for gradient projection method solving the full problem on the finest mesh.}
    \begin{tabular}{l|rr|rr|rr|rr|rr}
    \toprule
    level (vars) & \multicolumn{2}{c}{4\,(961)} & \multicolumn{2}{c}{5\,(3969)} & \multicolumn{2}{c}{6\,(16129)} & \multicolumn{2}{c}{7\,(65025)} & \multicolumn{2}{c}{8\,(261121)} \\\midrule
    \midrule
          smoother & \multicolumn{1}{c}{rate} & \multicolumn{1}{c}{feval} & \multicolumn{1}{c}{rate} & \multicolumn{1}{c}{feval} & \multicolumn{1}{c}{rate} & \multicolumn{1}{c}{feval} & \multicolumn{1}{c}{rate} & \multicolumn{1}{c}{feval} & \multicolumn{1}{c}{rate} & \multicolumn{1}{c}{feval} \\\midrule
    GSP   & 0.07 &       & 0.14 &       & 0.22 &       & 0.32 &       & 0.37 &  \\\midrule
    GP-1  & 0.18 & 71    & 0.33 & 107   & 0.50 & 180   & 0.80 & 410   & 0.86 & 711 \\
    GP-2  & 0.07 & 93    & 0.14 & 111   & 0.26 & 206   & 0.57 & 384   & 0.70 & 677 \\
    GP-3  & 0.03 & 92    & 0.08 & 142   & 0.17 & 239   & 0.35 & 387   & 0.55 & 806 \\
    GP-4  & 0.02 & 122   & 0.05 & 176   & 0.12 & 285   & 0.31 & 459   & 0.53 & 887 \\
    GP-5  & 0.01 & 160   & 0.03 & 211   & 0.08 & 306   & 0.23 & 488   & 0.34 & 912 \\\midrule
    GP only &    & 685   &      & 2361  &      & 9320  &      & 34133 &      & 127289 \\
    \bottomrule
    \end{tabular}%
  \label{tab:1}%
\end{table}%

\begin{figure}[h]
\begin{center}
 \resizebox{0.48\hsize}{!}
   {\includegraphics{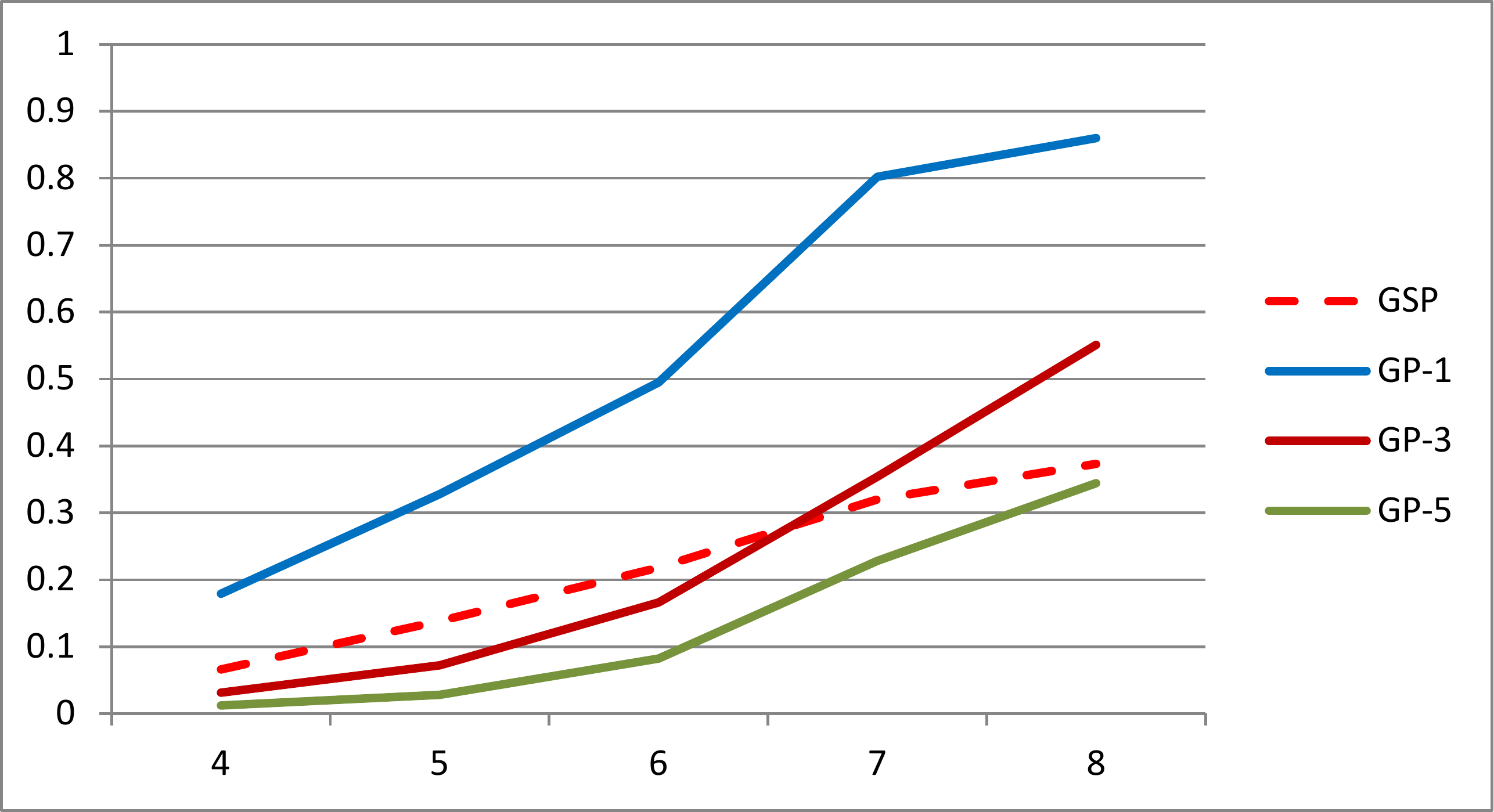}}\quad
 \resizebox{0.48\hsize}{!}
   {\includegraphics{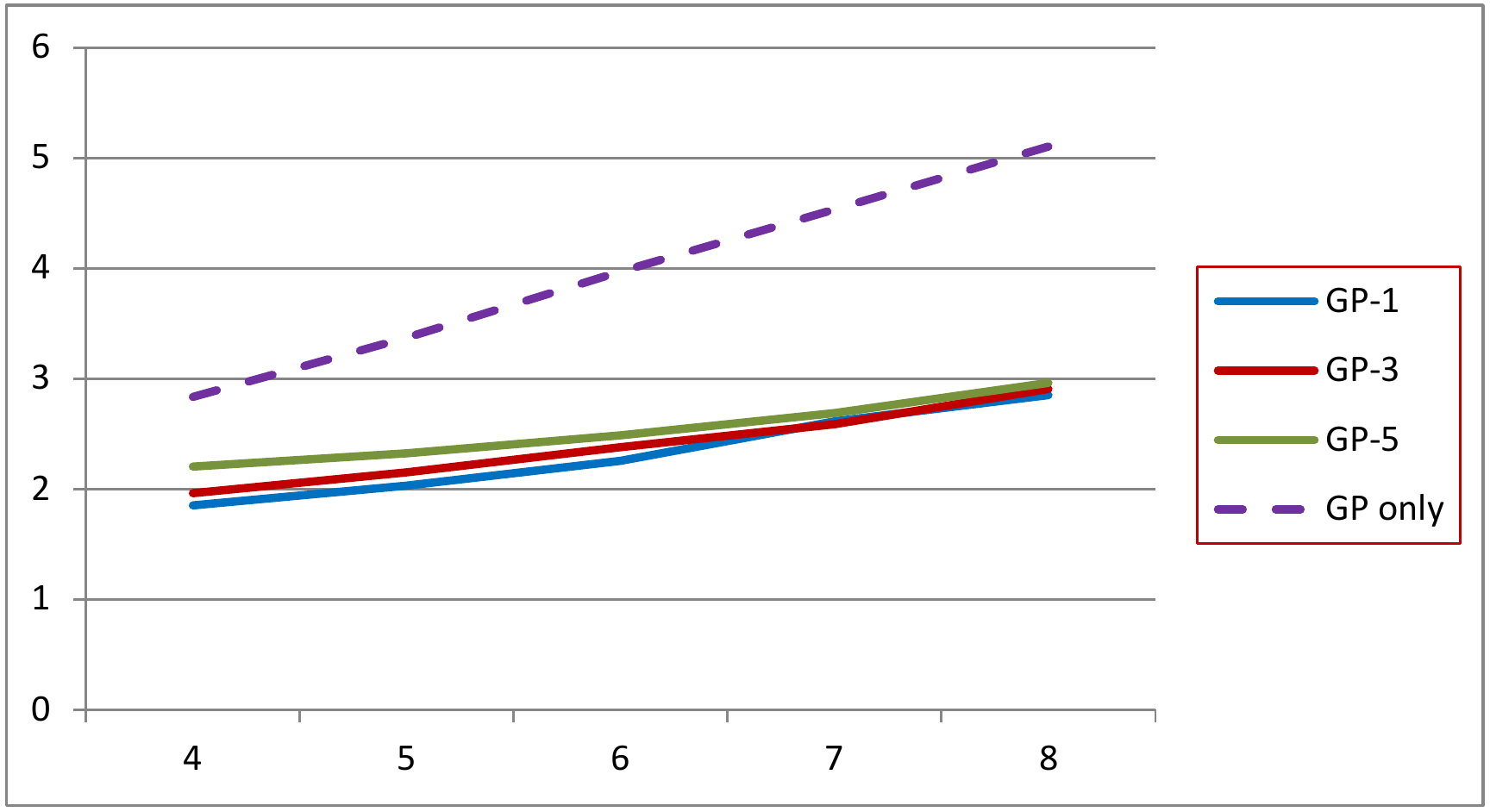}}
  \end{center}
  \caption{Example~\ref{ex:1}, rate of convergence (left) and function evaluations (right) for various smoothers as a function of the number of levels.}\label{fig:ex1a}
\end{figure}
The numbers from Table~\ref{tab:1} are graphically presented in Figure~\ref{fig:ex1a}: it shows the rate of convergence for the different smoothing steps, as it increases with the number of levels (left-hand figure). The right-hand figure presents the logarithm of function evaluations as a function of the number of levels. We can observe a rapid increase for the ``pure'' gradient projection method and a much smaller increase for the multigrid algorithm, almost independent of the number of smoothing steps.

\subsection{Example: non-quadratic obstacle problem}\label{ex:2}
Consider the following optimization problem in $\Omega\subset \RR^2$:
\begin{align*}
  &\min_{u\in H^1_0(\Omega)} {\cal J}(u):= \frac{1}{2}\int_\Omega \|\nabla u\|^2 - (ue^u-e^u)\;dx - \int_{\Omega} Fu\;dx\\
  &\mbox{subject to}\nonumber\\
  &\qquad\varphi \leq u \leq \psi,\quad \mbox{a.e.\ in}\ \Omega\,,
\end{align*}
with
$$
  \varphi(x_1,x_2) = -8 (x_1-7/16)^2-8(x_2-7/16)^2+0.2,\qquad \psi=0.5
$$
and
$$
  F(x_1,x_2)=\left(9 \pi^2+e^{(x_1^2-x_1^3)\sin(3\pi x_2)}(x_1^2-x_1^3)+6x_1-2\right) \sin(3\pi x_1)\,.
$$
\begin{figure}[h]
\begin{center}
 \resizebox{0.48\hsize}{!}
   {\includegraphics{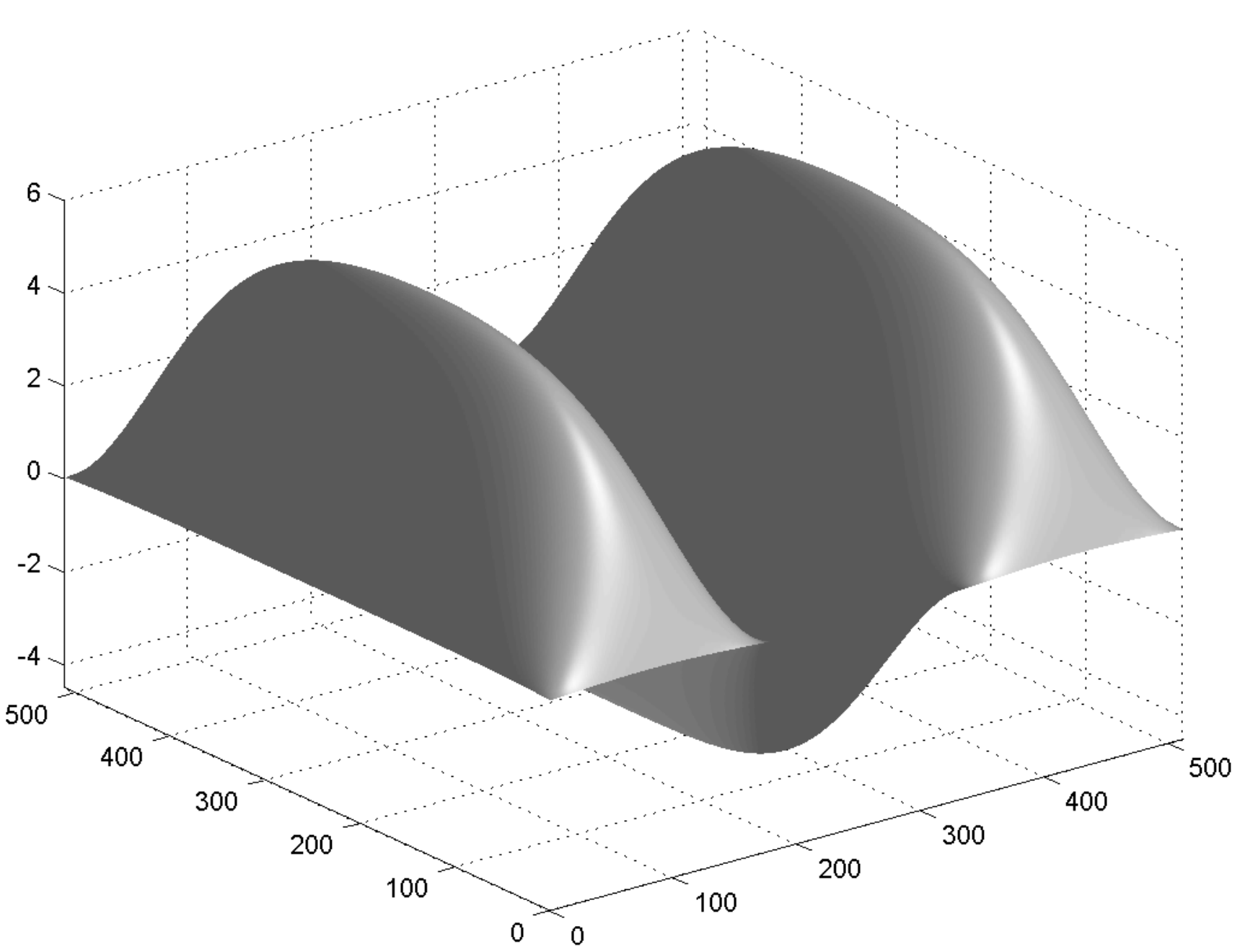}}\quad
 \resizebox{0.4\hsize}{!}
   {\includegraphics{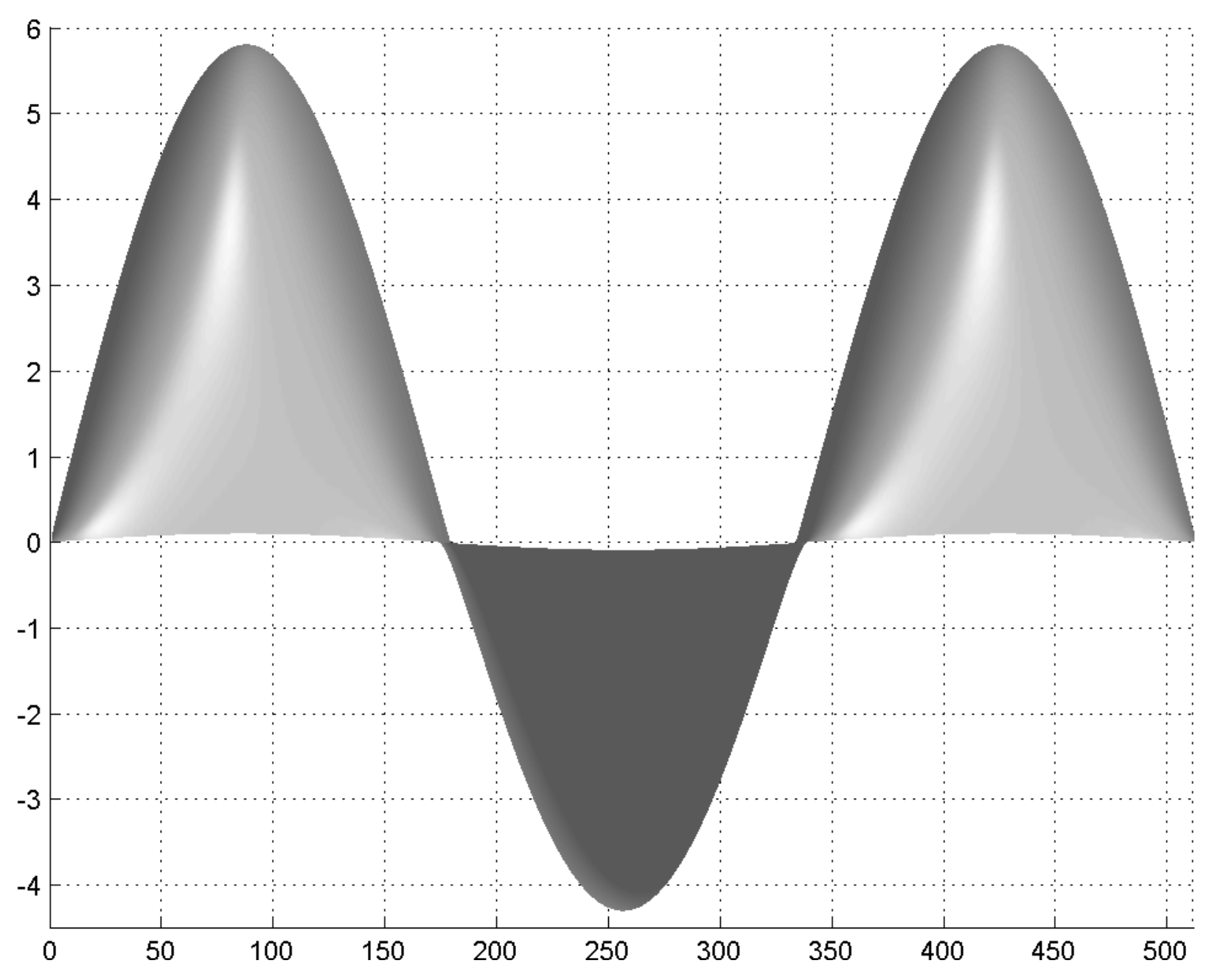}}\\
    \resizebox{0.48\hsize}{!}
   {\includegraphics{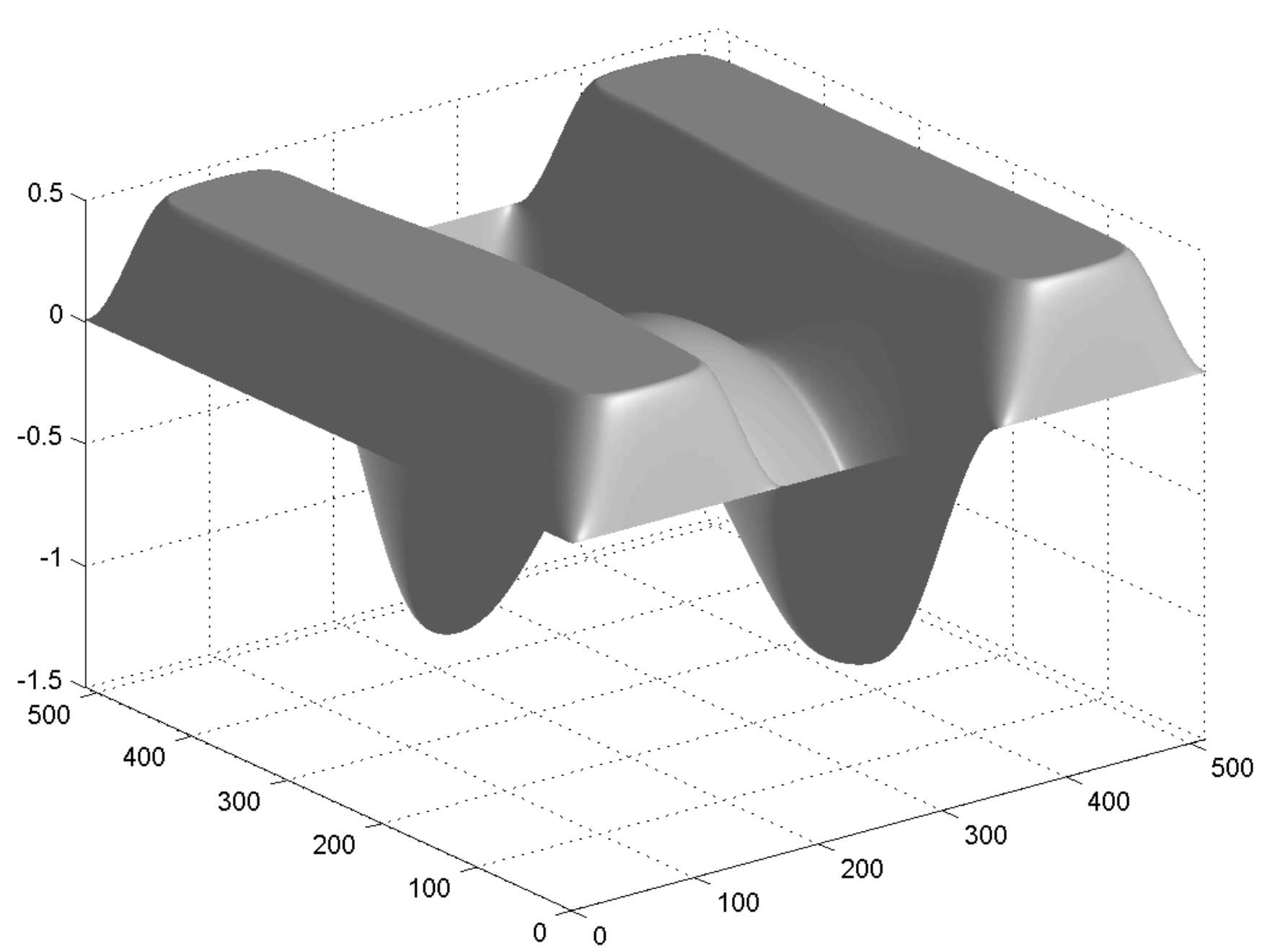}}\quad
 \resizebox{0.48\hsize}{!}
   {\includegraphics{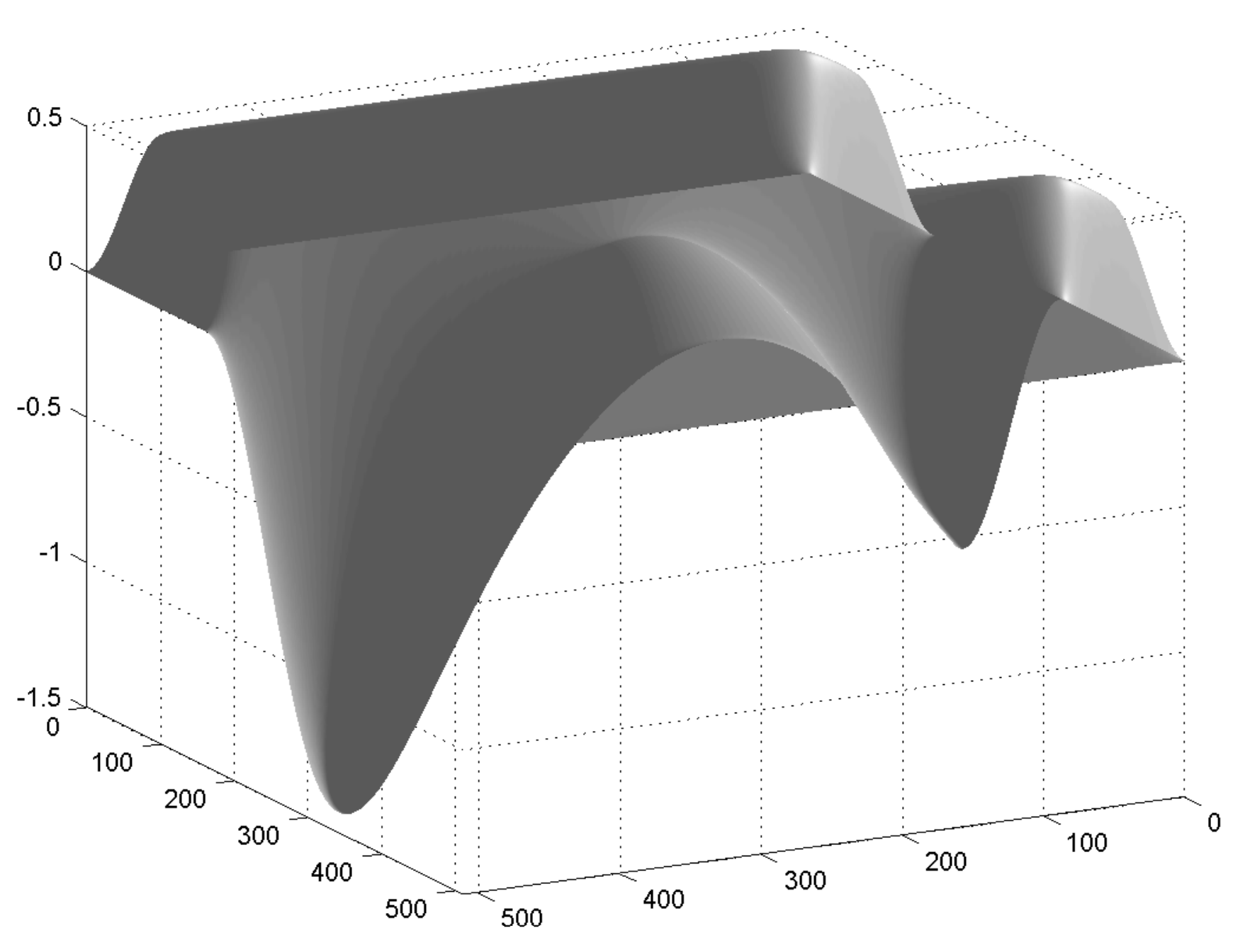}}
  \end{center}
  \caption{Example~\ref{ex:2}, solution of the unconstrained problem (top) and the constrained problem (bottom).}\label{fig:ex2}
\end{figure}
The unconstrained version of the problem is a nonlinear PDE studied in \cite[p.105]{briggs2000multigrid}.
Figure~\ref{fig:ex2} shows the solution of the unconstrained (top) and the constrained problem (bottom), both in two different views.

\begin{figure}[h]
\begin{center}
 \resizebox{0.48\hsize}{!}
   {\includegraphics{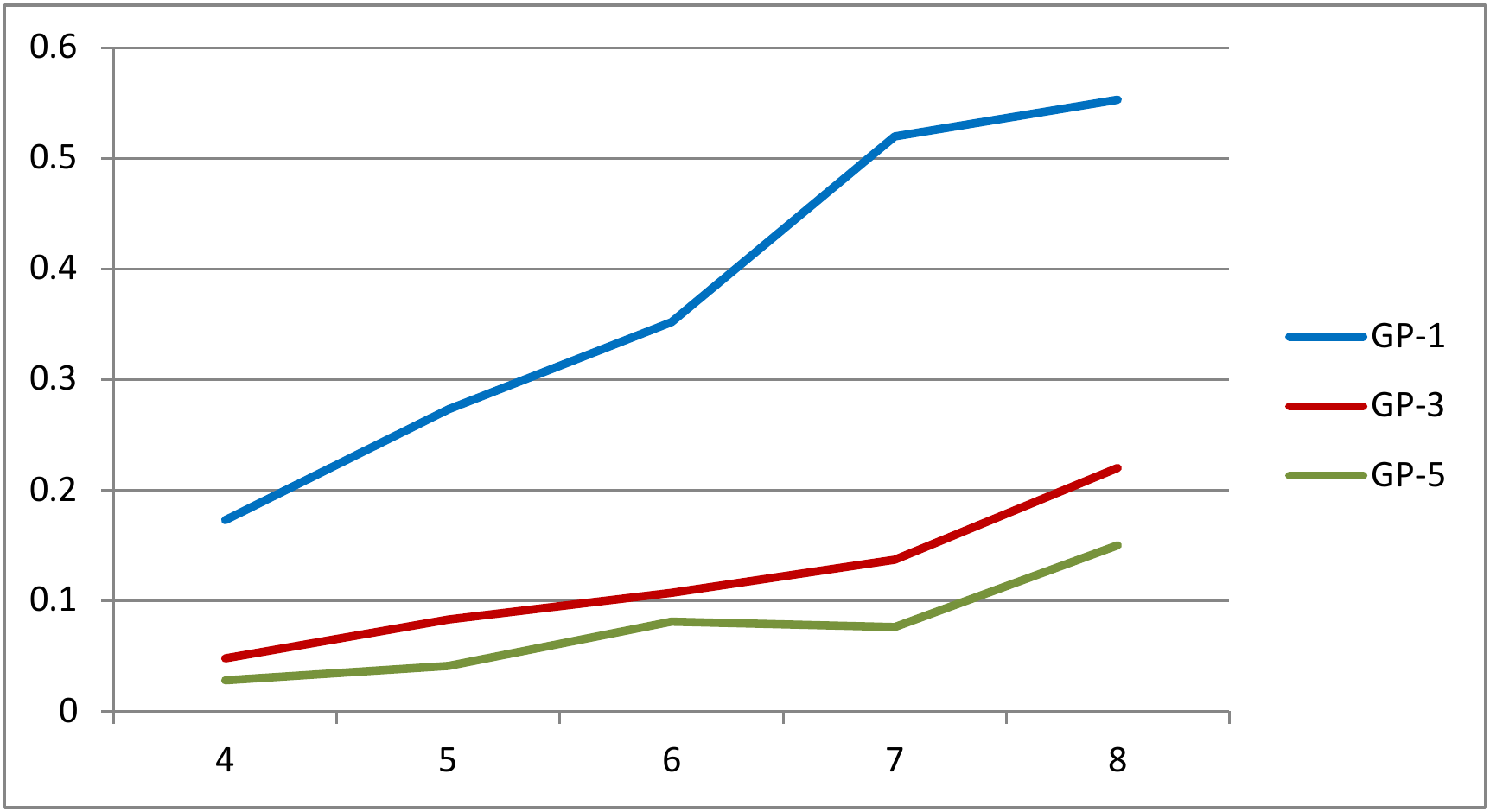}}\quad
 \resizebox{0.48\hsize}{!}
   {\includegraphics{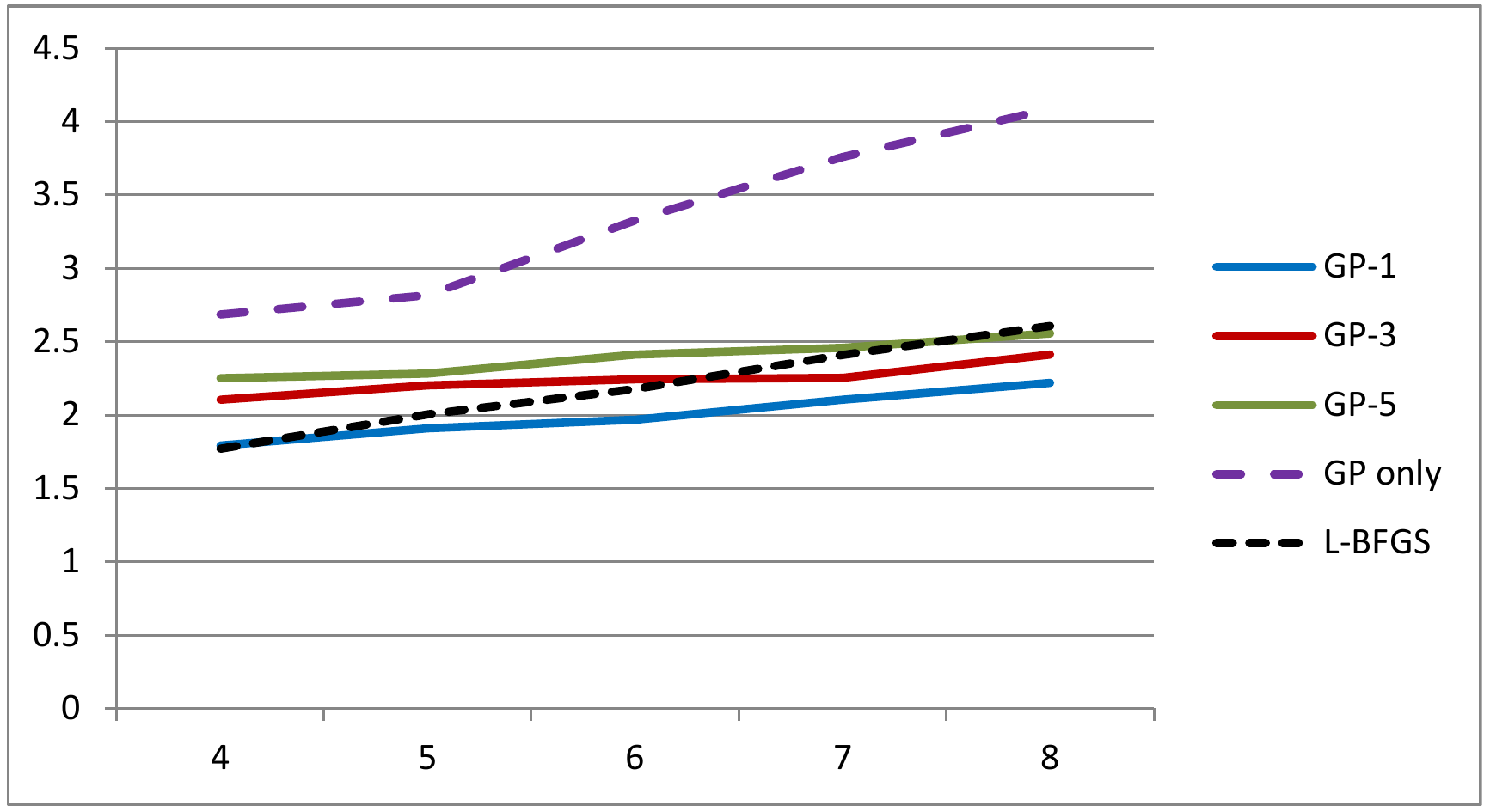}}\quad
  \end{center}
  \caption{Example~\ref{ex:2}, rate of convergence (left) and function evaluations (right) for various smoothers as a function of the number of levels.}\label{fig:ex2a}
\end{figure}

Just as in Example~\ref{ex:1}, Table~\ref{tab:2} together with Figure~\ref{fig:ex2a} present the results of the numerical experiments using Algorithm~2. We do not show a comparison of the GP smoother with the projected nonlinear Gauss-Seidel smoother, used, e.g., in \cite{hackbusch1983}. This is because the nonlinear GS algorithm needs the Hessian of the objective function and turns the algorithm into a second-order method. We can use a finite difference approximation of the Hessian but then the resulting code is extremely slow.
We do, however, compare the multigrid algorithm with one of the most efficient codes for bound-constrained nonlinear optimization, the L-BFGS-B code by Morales and Nocedal \cite{morales2011remark}. We can see that, with increasing size of the problem, the number of function evaluations grows faster in L-BFGS-B (though we should keep in mind that additional work needs to be done on coarse levels of the multigrid algorithm; see Remark~\ref{rem:1} below). And we should also keep in mind that, unlike in the multigrid algorithm, the function and gradient evaluation is not the only computationally expensive part of the L-BFGS-B code. We do not compare the CPU times, as L-BFGS-B is coded in Fortran.
\begin{table}[htbp]
  \centering
  \caption{Example~\ref{ex:2}, asymptotic rate of convergence and number of top-level function evaluations for 4--8 refinement levels. GP-$\nu$ stands for a $(\nu,\nu)$ V-cycle; ``GP only'' and ``L-BFGS-B'' for gradient projection method and the L-BFGS-B method, respectively, solving the full problem on the finest level.}
    \begin{tabular}{l|rr|rr|rr|rr|rr}
    \toprule
    level (vars) & \multicolumn{2}{c}{4\,(961)} & \multicolumn{2}{c}{5\,(3969)} & \multicolumn{2}{c}{6\,(16129)} & \multicolumn{2}{c}{7\,(65025)} & \multicolumn{2}{c}{8\,(261121)} \\\midrule
          smoother & \multicolumn{1}{c}{rate} & \multicolumn{1}{c}{feval} & \multicolumn{1}{c}{rate} & \multicolumn{1}{c}{feval} & \multicolumn{1}{c}{rate} & \multicolumn{1}{c}{feval} & \multicolumn{1}{c}{rate} & \multicolumn{1}{c}{feval} & \multicolumn{1}{c}{rate} & \multicolumn{1}{c}{feval} \\\midrule
    GP-1     & 0.17 & 62    & 0.27 & 81    & 0.35 & 93    & 0.52 & 127   & 0.55  & 166 \\
    GP-2     & 0.12 & 131   & 0.21 & 193   & 0.29 & 192   & 0.42 & 282   & 0.50  & 321 \\
    GP-3     & 0.05 & 127   & 0.08 & 159   & 0.11 & 175   & 0.14 & 179   & 0.22  & 258 \\
    GP-4     & 0.05 & 171   & 0.07 & 205   & 0.09 & 249   & 0.14 & 284   & 0.29  & 384 \\
    GP-5     & 0.03 & 178   & 0.04 & 192   & 0.08 & 259   & 0.08 & 288   & 0.15  & 360 \\\midrule
    GP only  &      & 485   &      & 656   &      & 2128  &      & 5746  &       & 12197 \\\midrule
    L-BFGS-B &      & 59    &      & 101   &      & 151   &      & 257   &       & 405 \\
    \bottomrule
    \end{tabular}%
  \label{tab:2}%
\end{table}%

\begin{remark}\label{rem:1}
To have a better idea about the amount of work required on the coarser levels we present CPU times for the largest problem ($\ell=9$) with GP-1. Using MATLAB's tic-toc commands, we measured the cumulative times spent in the iterative method (the smoother) on every level; the times are given in seconds. The time spent on the finest mesh was 6.240, while the time spent on all other meshes was 1.824. The latter number is a sum of $(0.231, 0.036, 0.029, 0.031, 0.043, 0.108, 0.250, 1.097)$ corresponding to the coarsest up to the second finest mesh, respectively. Recall that the problem on the coarsest mesh is solved to high accuracy.
\end{remark}

\subsection{Example: minimal surface problem}\label{ex:3}
\begin{figure}[h]
\begin{center}
 \resizebox{0.48\hsize}{!}
   {\includegraphics{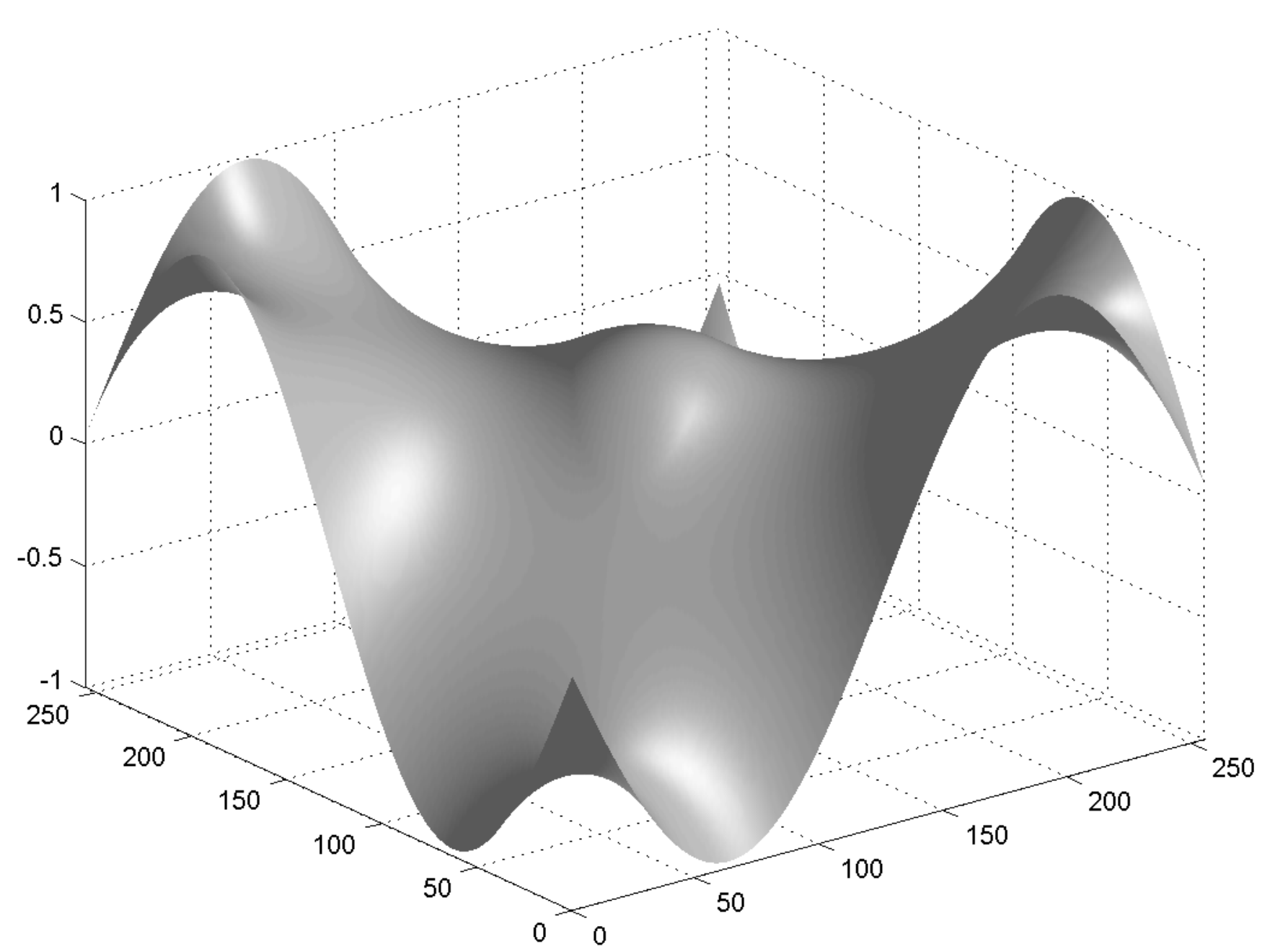}}\quad
 \resizebox{0.4\hsize}{!}
   {\includegraphics{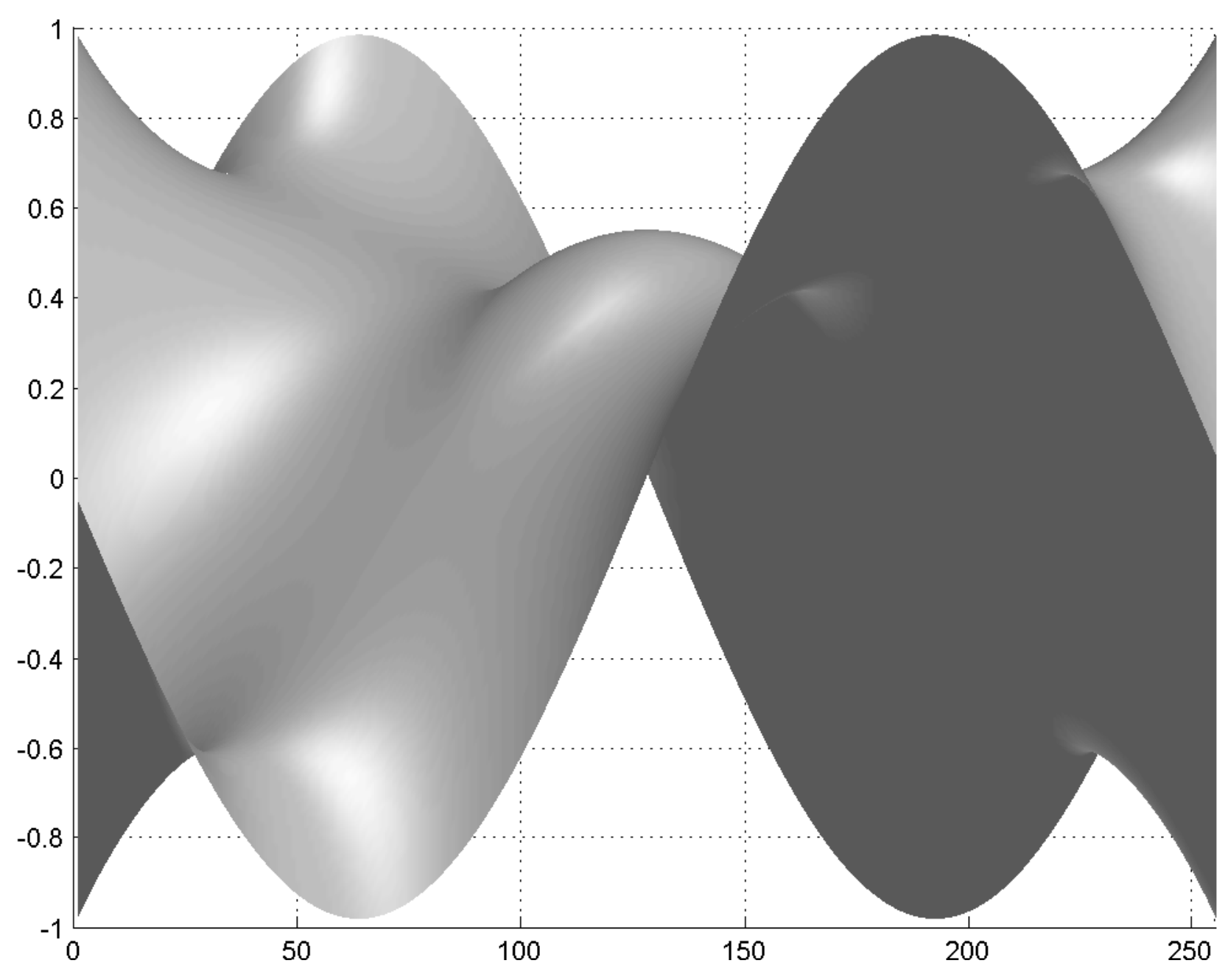}}
  \end{center}
  \caption{Example~\ref{ex:3}, solution.}\label{fig:ex4}
\end{figure}

Our next example is the minimal surface problem
\begin{align*}
  &\min_{u\in H^1(\Omega)} {\cal J}(u):= \int_\Omega \sqrt{1+\|\nabla u\|^2}\;dx \\
  &\mbox{subject to}\nonumber\\
  &\qquad u(x_1,x_2) = u_\Gamma(x_1,x_2)\ \mbox{for}\ (x_1,x_2)\in \partial\Omega\\
  &\qquad \varphi \leq u \leq \psi,\quad \mbox{a.e.\ in}\ \Omega\,,
\end{align*}
with the boundary function (see \cite{gratton2010numerical})
$$
  u_\Gamma(x_1,0)=\omega,\ u_\Gamma(1,x_2)=-\omega,\ u_\Gamma(x_1,1)=-\omega,\ u_\Gamma(0,x_2)=\omega,\ \omega=-\sin(2\pi \xi)
$$
and the parabolic lower bound
$$
  \varphi(x_1,x_2) = -8(x_1-0.5)^2-8(x_2-0.5)^2+0.55\,.
$$
The upper bound function $\psi$ is set to infinity. The solution is shown in two different views in Figure~\ref{fig:ex4}.

\begin{figure}[h]
\begin{center}
 \resizebox{0.48\hsize}{!}
   {\includegraphics{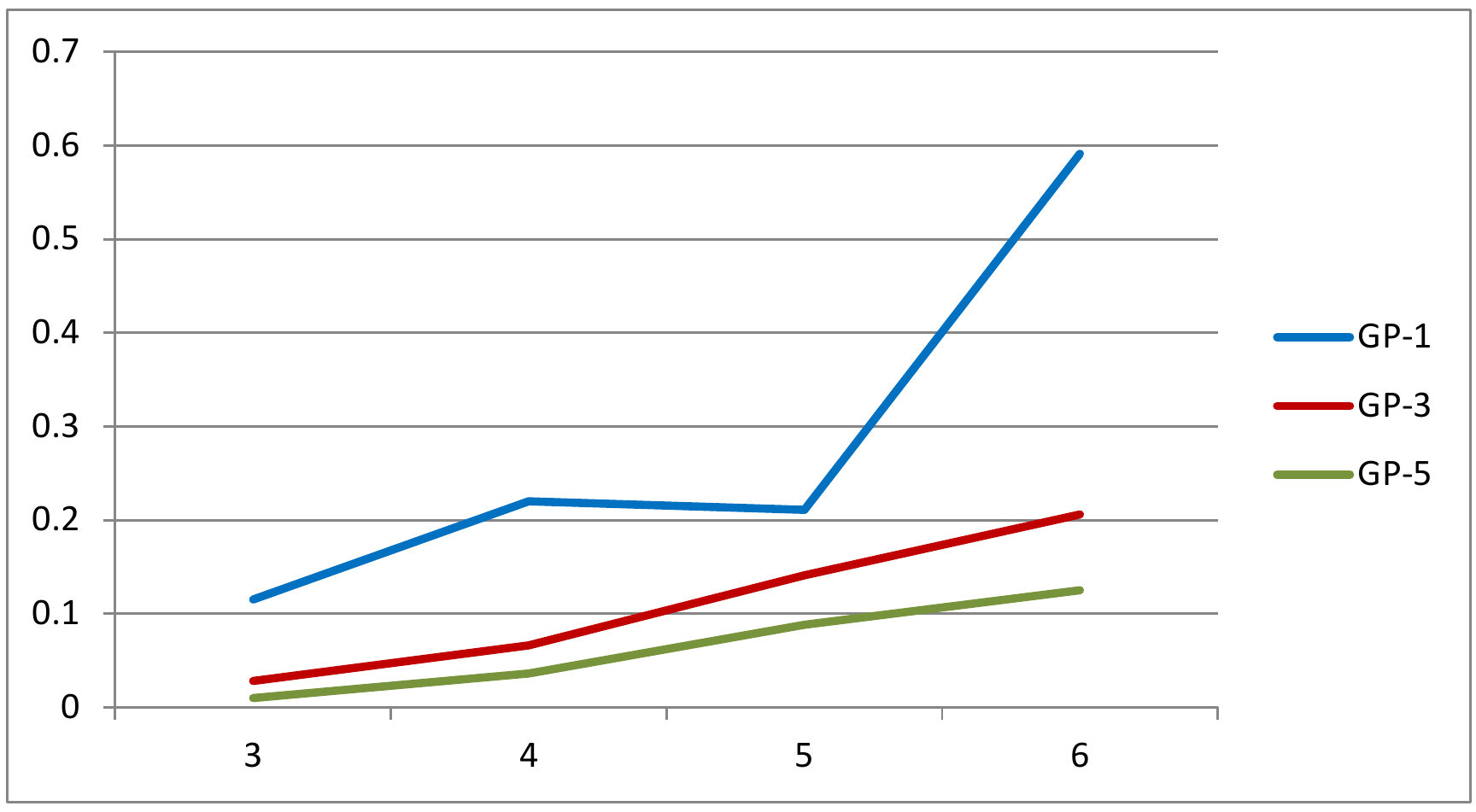}}\quad
 \resizebox{0.48\hsize}{!}
   {\includegraphics{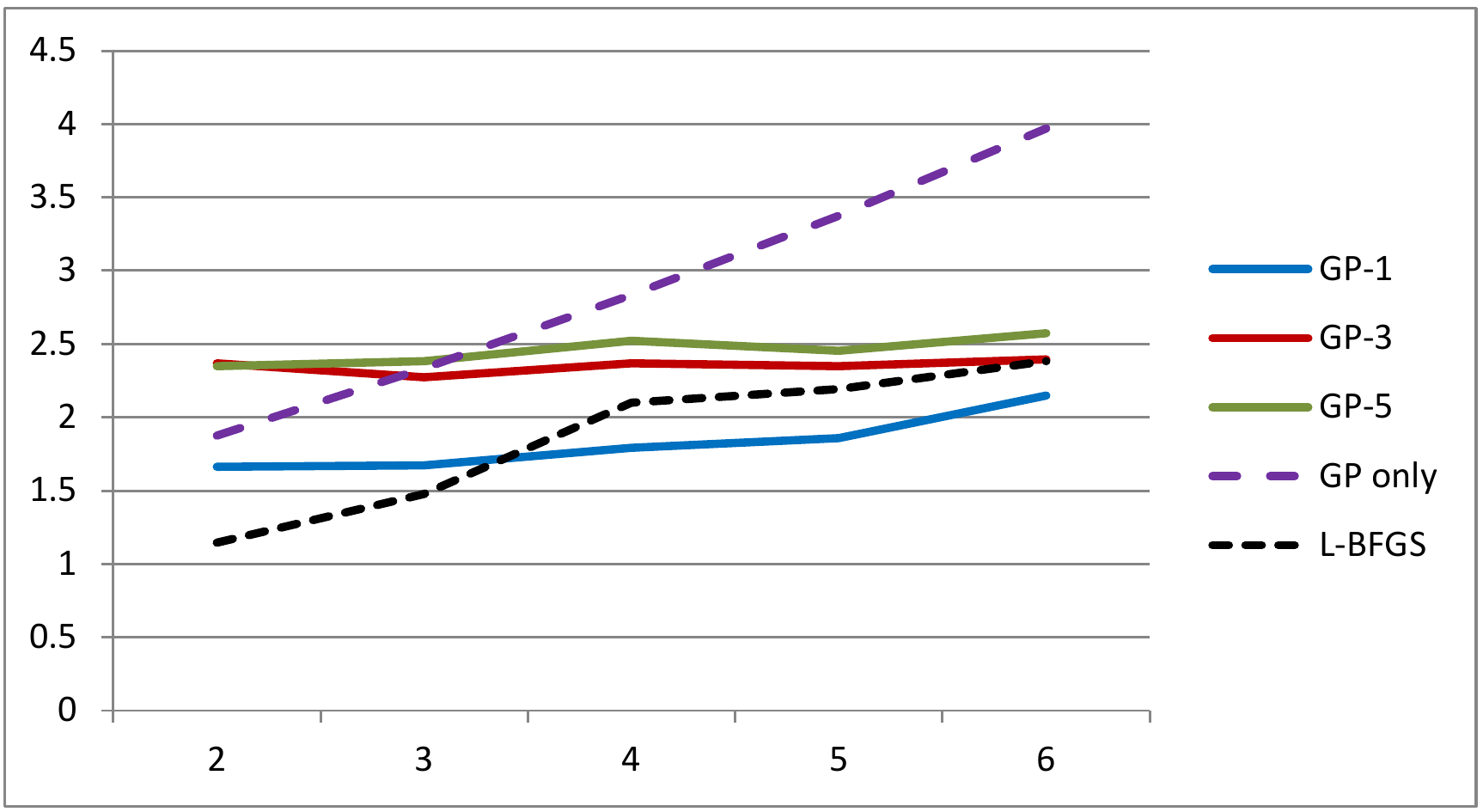}}
  \end{center}
  \caption{Example~\ref{ex:3}, rate of convergence (left) and function evaluations (right) for various smoothers as a function of the number of levels.}\label{fig:ex3a}
\end{figure}

As before, Table~\ref{tab:3} together with Figure~\ref{fig:ex3a} present the results of the numerical experiments using Algorithm~2. Notice that the function and gradient evaluation for this problem is much more expensive than in the other examples and we were not able to obtain the exact solution to the finest problems in a reasonable time by L-BFGS-B. That is why we only present results for levels 3--7.
\begin{table}[htbp]
  \centering
  \caption{Example~\ref{ex:3}, asymptotic rate of convergence and number of top-level function evaluations for 2--6 refinement levels. GP-$\nu$ stands for a $(\nu,\nu)$ V-cycle; ``GP only'' and ``L-BFGS-B'' for gradient projection method and the L-BFGS-B method, respectively, solving the full problem on the finest level.}
    \begin{tabular}{l|rr|rr|rr|rr|rr}
    \toprule
    level (vars) & \multicolumn{2}{c}{2\,(81)} & \multicolumn{2}{c}{3\,(289)} & \multicolumn{2}{c}{4\,(1089)} & \multicolumn{2}{c}{5\,(4225)} & \multicolumn{2}{c}{6\,(16641)}\\
    \midrule
          smoother & \multicolumn{1}{c}{rate} & \multicolumn{1}{c}{feval} & \multicolumn{1}{c}{rate} & \multicolumn{1}{c}{feval} & \multicolumn{1}{c}{rate} & \multicolumn{1}{c}{feval} & \multicolumn{1}{c}{rate} & \multicolumn{1}{c}{feval} & \multicolumn{1}{c}{rate} & \multicolumn{1}{c}{feval} \\\midrule
    GP-1     & 0.118 & 46    & 0.115 & 47    & 0.22 & 62    & 0.21 & 72    & 0.60 & 141 \\
    GP-2     & 0.079 & 125   & 0.055 & 157   & 0.13 & 176   & 0.18 & 111   & 0.43 & 171 \\
    GP-3     & 0.029 & 234   & 0.028 & 188   & 0.07 & 234   & 0.14 & 223   & 0.21 & 249 \\
    GP-4     & 0.011 & 156   & 0.021 & 205   & 0.05 & 165   & 0.06 & 170   & 0.19 & 276 \\
    GP-5     & 0.004 & 223   & 0.010 & 241   & 0.04 & 332   & 0.09 & 284   & 0.13 & 374 \\\midrule
    GP only  &       & 75    &       & 216   &      & 685   &      & 2361  &      & 9320 \\\midrule
    L-BFGS-B &       & 14    &       & 30    &      & 126   &      & 156   &      & 242 \\
    \bottomrule
    \end{tabular}%
  \label{tab:3}%
\end{table}%

\subsection{Example: obstacle problem with an equality constraint}\label{ex:4}
\begin{figure}[h]
\begin{center}
 \resizebox{0.45\hsize}{!}
   {\includegraphics{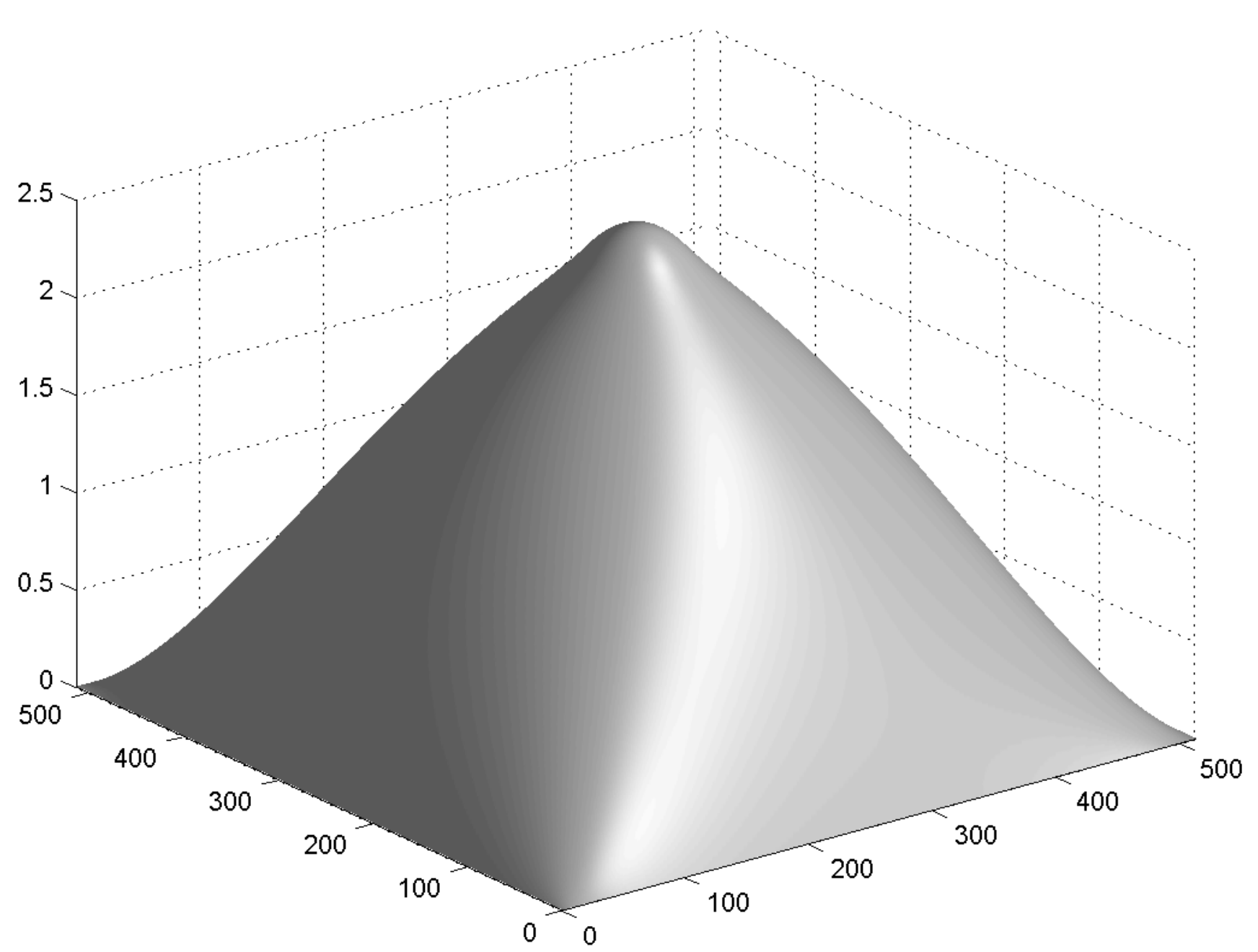}}\quad
 \resizebox{0.45\hsize}{!}
   {\includegraphics{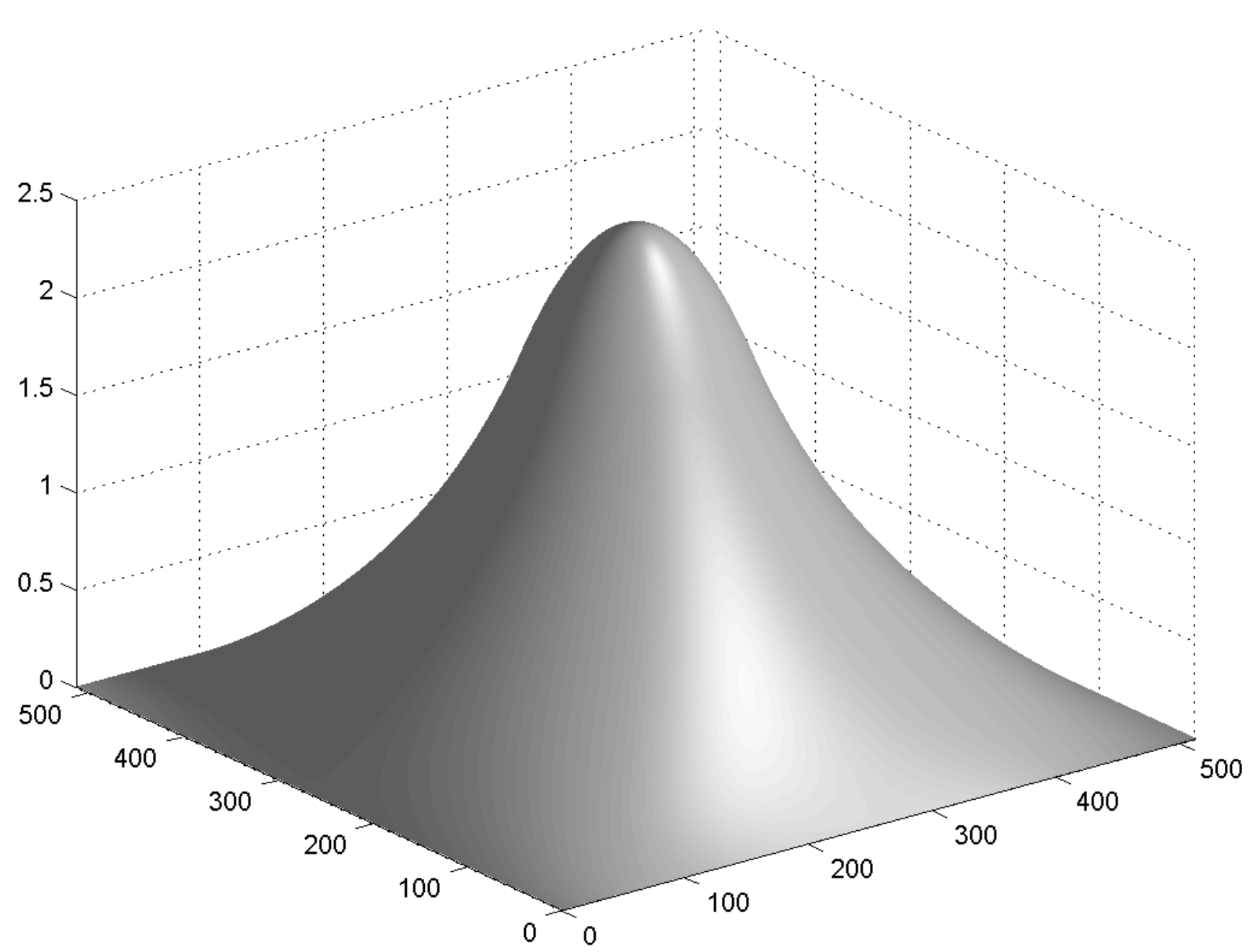}}
  \end{center}
  \caption{Example~\ref{ex:4}, solution with (left) and without (right) the equality constraint.}\label{fig:ex5}
\end{figure}

Finally, let us consider an example with an obstacle and an additional equality constraint. The problem stems from the nonlinear PDE
\begin{alignat*}{2}
  -\triangle u - u^2 & = f(x)\quad &\mbox{in}&\ \Omega\\
  u&=0\quad &\mbox{on}&\ \partial\Omega
\end{alignat*}
and can be formulated as the following optimization problem
\begin{align*}
  &\min_{u\in H^1_0(\Omega)} {\cal J}(u):= \frac{1}{2}\int_\Omega \left(\|\nabla u\|^2 - \frac{1}{3}u^3\right)\;dx - \int_{\Omega} Fu\;dx\\
  &\mbox{subject to}\nonumber\\
  &\qquad \varphi \leq u \quad \mbox{a.e.\ in}\ \Omega\\
  &\qquad \int_\Omega u\;dx = 1\,,
\end{align*}
with $F\equiv 0$ and
$$
  \varphi(x_1,x_2) = -32 (x_1-0.5)^2-32(x_2-0.5)^2+2.5\,.
$$
Figure~\ref{fig:ex5} (left) shows the solution and a comparison with the solution of the same problem without the equality constraint (right). In the unconstrained case, the optimal solution gives $\int_\Omega u\;dx = 0.62$. So, in order to satisfy the equality constraint, the unconstrained solution has been ``inflated''.

\begin{figure}[h]
\begin{center}
 \resizebox{0.48\hsize}{!}
   {\includegraphics{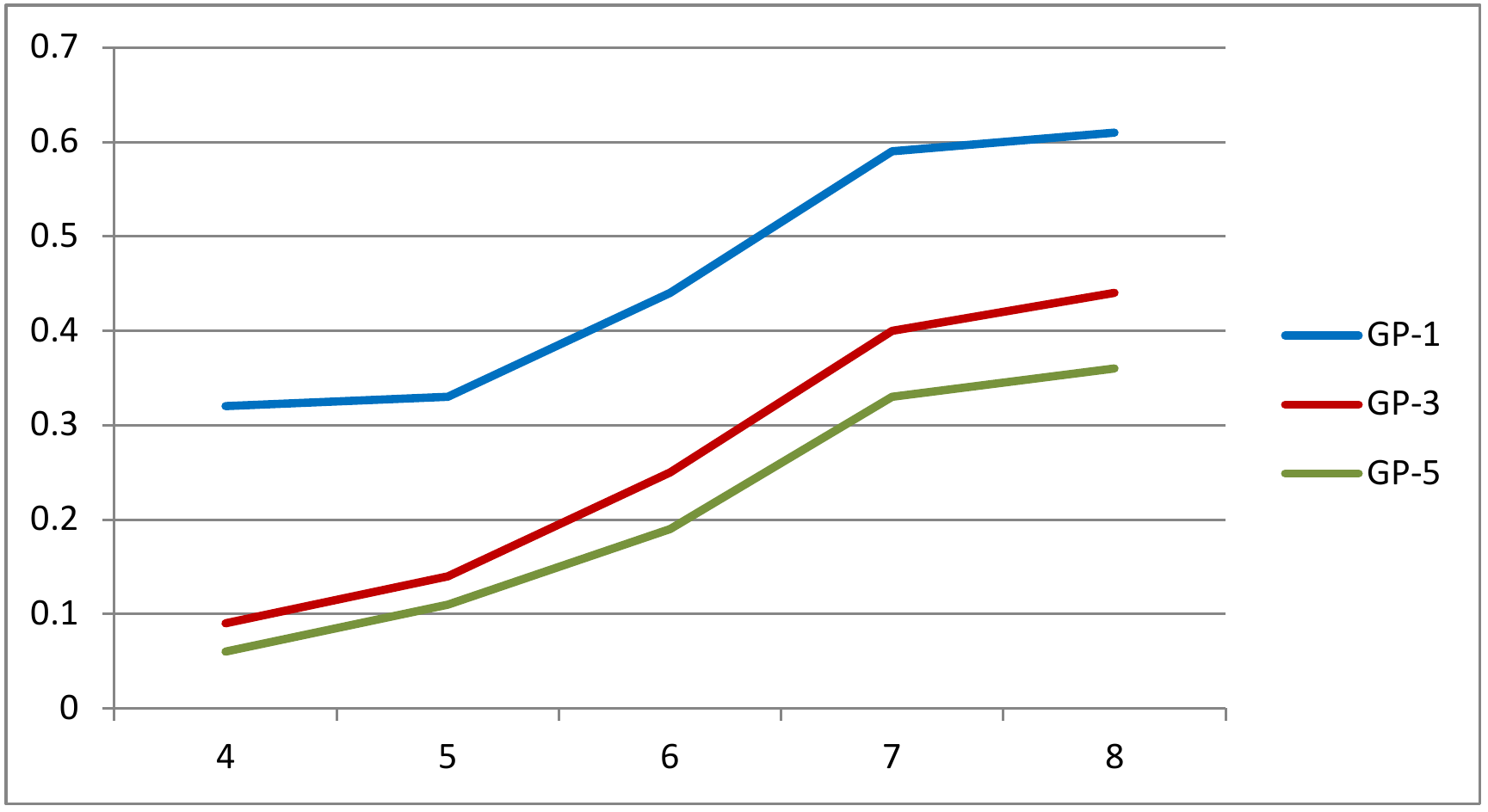}}\quad
 \resizebox{0.48\hsize}{!}
   {\includegraphics{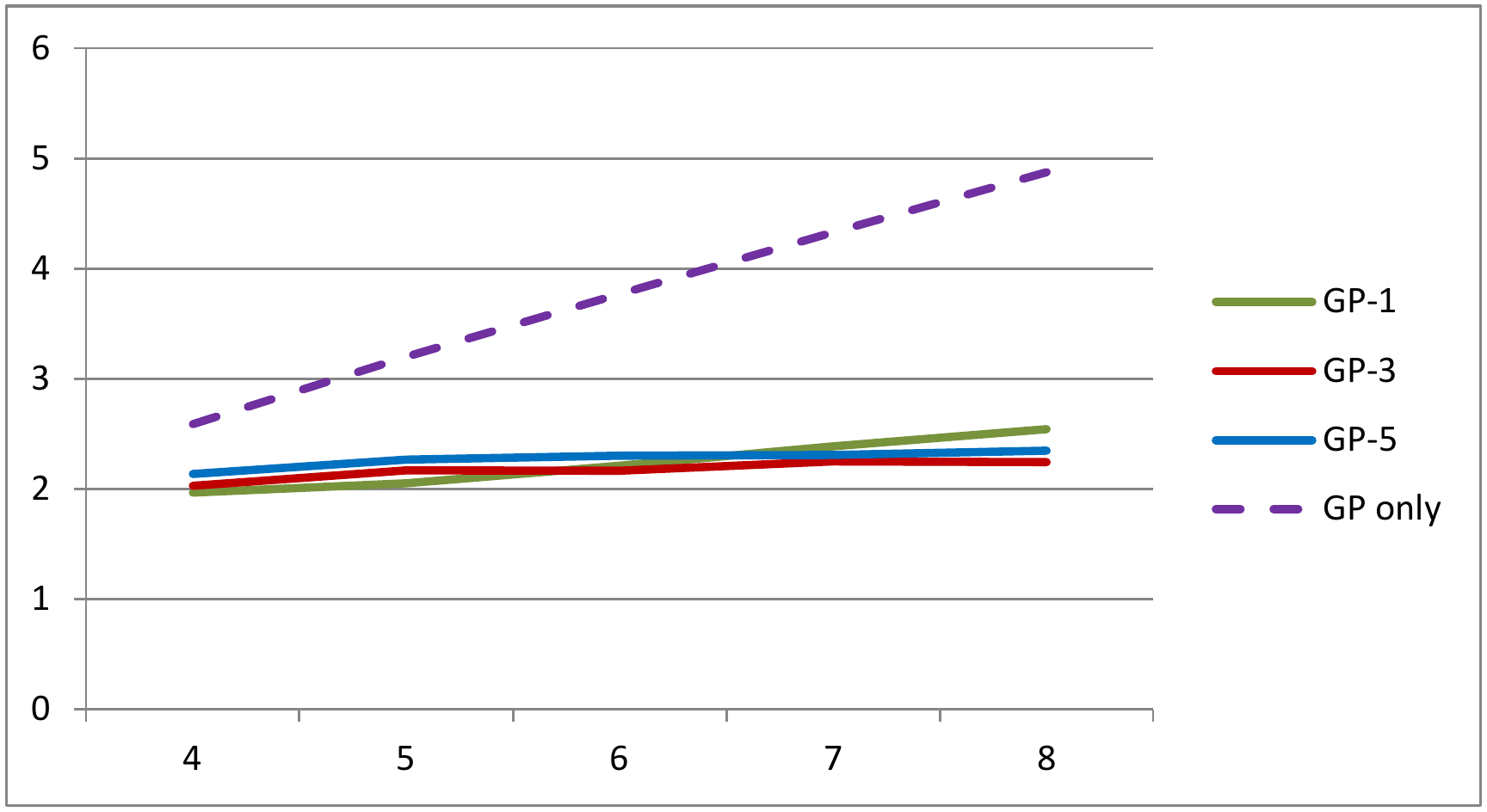}}
  \end{center}
  \caption{Example~\ref{ex:4}, rate of convergence (left) and function evaluations (right) for various smoothers as a function of the number of levels.}\label{fig:ex4a}
\end{figure}

Table~\ref{tab:4} together with Figure~\ref{fig:ex4a} present the results of the numerical experiments using Algorithm~3 with the additional handling of the equality constraint (Section~\ref{sec:eq}). The explanation is the same as in the previous examples.
\begin{table}[htbp]
  \centering
  \caption{Example~\ref{ex:4}, asymptotic rate of convergence and number of top-level function evaluations for 4--8 refinement levels. GP-$\nu$ stands for a $(\nu,\nu)$ V-cycle; ``GP only'' for gradient projection method solving the full problem on the finest level.}
    \begin{tabular}{l|rr|rr|rr|rr|rr}
    \toprule
    level (vars) & \multicolumn{2}{c}{4\,(961)} & \multicolumn{2}{c}{5\,(3969)} & \multicolumn{2}{c}{6\,(16129)} & \multicolumn{2}{c}{7\,(65025)} & \multicolumn{2}{c}{8\,(261121)} \\
    \midrule
         smoother  & \multicolumn{1}{c}{rate} & \multicolumn{1}{c}{feval} & \multicolumn{1}{c}{rate} & \multicolumn{1}{c}{feval} & \multicolumn{1}{c}{rate} & \multicolumn{1}{c}{feval} & \multicolumn{1}{c}{rate} & \multicolumn{1}{c}{feval} & \multicolumn{1}{c}{rate} & \multicolumn{1}{c}{feval} \\ \midrule
    GP-1    & 0.32  & 93    & 0.33  & 113   & 0.44  & 163   & 0.59  & 244   & 0.61  & 350 \\
    GP-2    & 0.11  & 88    & 0.25  & 120   & 0.29  & 129   & 0.51  & 183   & 0.54  & 182 \\
    GP-3    & 0.09  & 107   & 0.14  & 148   & 0.25  & 147   & 0.4   & 178   & 0.44  & 176 \\
    GP-4    & 0.07  & 153   & 0.14  & 182   & 0.23  & 186   & 0.36  & 224   & 0.44  & 191 \\
    GP-5    & 0.06  & 137   & 0.11  & 185   & 0.19  & 202   & 0.33  & 204   & 0.36  & 223 \\\midrule
    GP only &       & 388   &       & 1586  &       & 5907  &       & 21372 &       & 75258 \\
    \bottomrule
    \end{tabular}%
  \label{tab:4}%
\end{table}%

\begin{remark}
Notice that, in the presence of the equality constraint, we can no longer use Algorithm 5 as a smoother, as the gradient-based line search would not lead to a convergent algorithm. Instead, we use a standard projected gradient method with backtracking Armijo line search. To find the projection on the feasible set, we now have to solve a convex quadratic programming problem. Moreover, this problem has to be solved to a high precision, because we need to identify the active constraints in (\ref{eq:5})--(\ref{eq:6}). In our implementation, we have used the Gurobi solver for this purpose \cite{gurobi}.
\end{remark}

\subsection{To truncate or not to truncate}\label{sec:truncate}
\begin{figure}[h]
\begin{center}
 \resizebox{0.48\hsize}{!}
   {\includegraphics{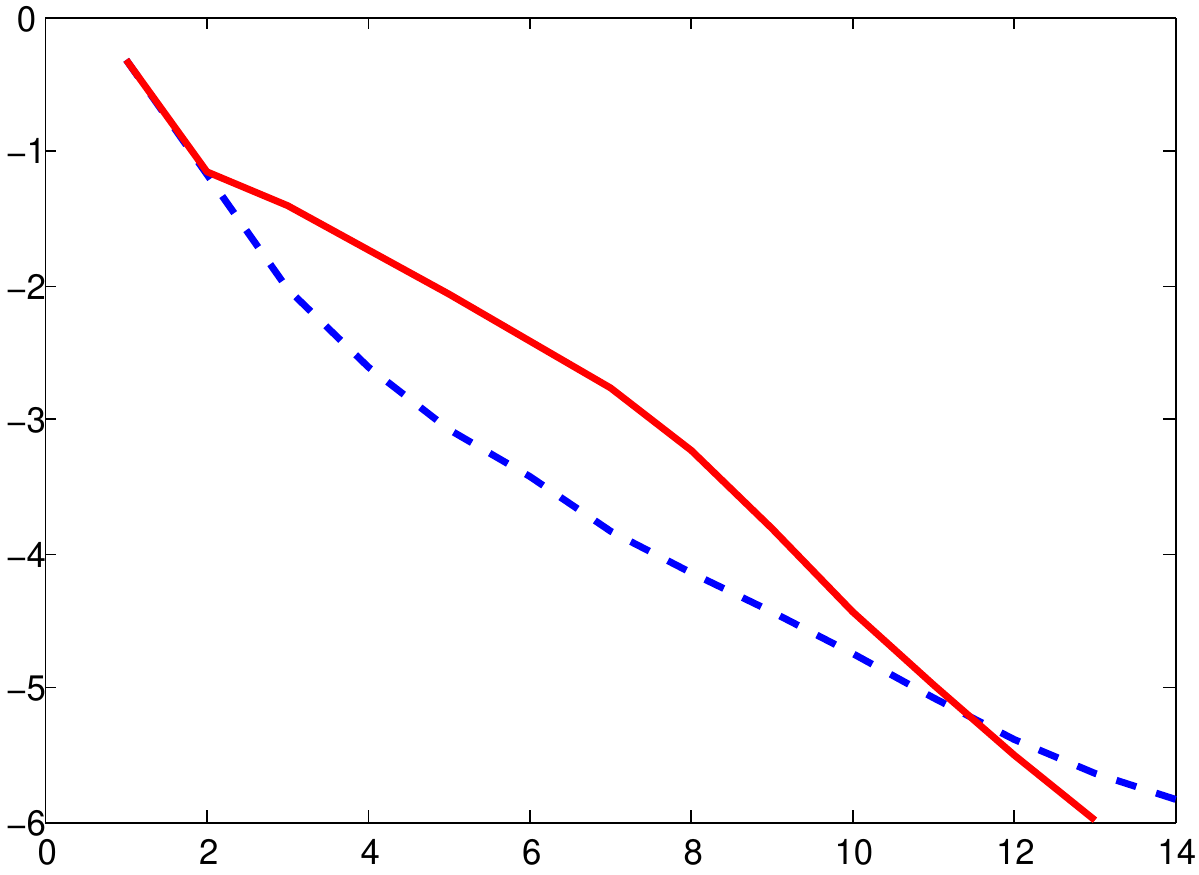}}\quad
  \end{center}
  \caption{Example~\ref{ex:1}, eight refinement levels. Convergence curves (iterations vs logarithm of the error) for Algorithm~2 with truncation (full line) and Algorithm~3 without truncation (dashed line).}\label{fig:trunc}
\end{figure}

The tables in \cite{Graser2009} (and partly in the previous section) show the clear advantage of truncation: the higher asymptotic rate of convergence as compared to Algorithm~3 without truncation. However, a typical user may not be interested in asymptotic rate but in fast convergence in the first iterations. And here Algorithm~3 can be the winner. Figure~\ref{fig:trunc} presents the convergence curves for Example~\ref{ex:1} with 8 refinement levels and $\nu=5$. The dashed line is for Algorithm~3 (no truncation) while the full line for Algorithm~2 (truncation). We can see a typical behaviour of the truncated algorithm: it starts slowly, tries to find the exact active set and, once this is found, the algorithm speeds up (for more details, see \cite{Graser2009}). However, the total amount of work (represented in this case by the integral of the convergence curve) to reach the required accuracy is actually lower for the asymptotically slower algorithm without truncation.

\section{Conclusions}
We have presented a version of the multigrid method for convex optimization problems with bound constraints and a possible single linear equality constraint. The method only needs gradient information, unlike similar published algorithms. We have shown that the projected gradient method can serve well as a smoother and that only a very small number of pre- and post-smoothing iterations is needed to obtain an efficient algorithm. The main advantage of the proposed method is thus in its low computational complexity and low memory requirements.

As an interesting by-product for unconstrained problems, we obtained a first-order method able to solve large scale problems efficiently and \emph{to high accuracy}, which is rather untypical in today's realm of first-order methods designed to solve very large scale convex problems though only to some 2--3 digits of accuracy.

The natural question arises about more general constraints. The authors devoted significant effort to the generalization of this method to the topology optimization problem \cite{bendsoe-sigmund}, a convex problem with bound constraints and a single linear equality constraint.  This problem, however, works with two kinds of variables that need to be discretized and prolonged/restricted on different refinement levels. It turns out that the interplay of these two kinds of discretization brings major technical difficulties. Also, it seems that more general constraints may increase the complexity of the formulas for constraint restriction. So at the time of writing this article, we think that more complex constraints could be better handled by traditional optimization methods (SQP or interior point) and multigrid could then be used for the solution of resulting linear systems.

Nevertheless, the class of unconstrained and bound-constrained convex optimization problems is very large and we believe that the presented method, whenever applicable, is one of the most efficient approaches to their solution.

\section*{Acknowledgements}
We would like to thank three anonymous reviewers for their valuable comments.
This work has been partly supported by Iraqi Ministry of Higher Education and Scientific Research,
Republic of Iraq, by the EU FP7 project AMAZE, and by the Grant Agency of the Czech Republic through project GAP201-12-0671. Their support is greatly acknowledged.

\bibliographystyle{gOMS}
\bibliography{mgm_paper}

\end{document}